\documentclass[12pt]{article}
\usepackage{amssymb,amsbsy}
\usepackage{color}
\begin{document}
\begin{center}
{\Large\bf  Constrained discounted stochastic games}
 \end{center}

\begin{center}
{\bf \large  Anna Ja\'skiewicz$^{a}$, Andrzej S. Nowak$^{b}$ }  
\end{center}
\noindent$^{a}$Faculty of Pure and Applied Mathematics,  Wroc{\l}aw University of Science and Technology, 
Wroc{\l}aw, Poland, 
{\footnotesize {\it email: anna.jaskiewicz@pwr.edu.pl}}\\
\noindent $^{b}$Faculty of Mathematics, Computer Science and Econometrics, University of Zielona G\'ora,
Zielona G\'ora, Poland,
{\footnotesize {\it email:  a.nowak@wmie.uz.zgora.pl}}\\

\begin{center}
\today
\end{center}
\vspace{3mm}
\noindent{\bf Abstract.} 
In this paper, we consider a large class of constrained non-cooperative stochastic Markov 
games with countable state spaces and discounted cost criteria. 
In one-player case, i.e., constrained discounted Markov decision  models,
it is possible to formulate a static optimisation problem 
whose solution determines a stationary optimal 
strategy ({\it alias} control or policy)  
in the dynamical infinite horizon model. 
This solution lies in the compact convex set  of 
all occupation measures induced by  strategies,  
defined on the set of state-action pairs.   
In case of $n$-person   discounted games the occupation measures 
are induced by strategies of all players. 
Therefore, it is difficult to generalise  the approach for constrained discounted
Markov decision processes  directly. It is not clear how to define the domain 
for the best-response correspondence whose fixed point induces a stationary
equilibrium in the Markov game. This domain should be the Cartesian 
product of compact convex sets in 
locally convex topological vector spaces.  
One of our main results shows how to overcome this difficulty
and define a constrained non-cooperative static 
game whose Nash equilibrium induces by a stationary
Nash equilibrium in the Markov game. This is done for games with bounded cost 
functions and positive initial state distribution.
An extension to a class of Markov games with unbounded costs and 
arbitrary initial state distribution relies on approximation of the unbounded 
game by bounded ones with positive initial state distributions.
In the unbounded case, we assume the uniform  integrability
of the discounted costs with respect to all probability measures induced 
by strategies of the players, defined on the space of plays (histories) of the game.
Our assumptions are weaker than those applied in 
earlier  works on discounted dynamic programming or stochastic games
using so-called weighted norm approaches. \\


 
\section{Introduction}

Constrained Markov decision processes arise in situations, in which a controller has  many objectives. For example,    
when she or he  wants to minimise one type of cost while keeping other
costs lower than some given bounds. 
Such situations appear very often in computer networks and data communications. 
The selected applications the reader may find in  the papers of 
Lazar \cite{laz},   Hordijk and Spieksma \cite{hs},  Ross and
Chen \cite{rc}  or  Feinberg and Reiman \cite{fr}. The theory of constrained
Markov decision processes  goes back to Derman and Klein \cite{dk}. It was further developed for finite state space models by  
Kallenberg \cite{kal}.
For the literature dealing with the discounted costs (or rewards), 
the reader is referred to the works of Altman \cite{alt}, Borkar \cite{bor}, 
Feinberg et al. \cite{auto}, Feinberg and Shwartz \cite{fs}, Sennott \cite{sen}
and the books by Altman \cite{a} and Piunovskiy  \cite{piu}.

Unconstrained non-cooperative $n$-person discounted stochastic Markov games with finite state spaces
were first studied by Fink \cite{fink}, Takahashi \cite{taka} and Sobel \cite{sobel}. Their results were extended to 
countable state space games by Federgruen \cite{f}. The proofs of the existence of stationary Nash equilibria   
in unconstrained discounted Markov games  are based on the Kakutani-Fan-Glicksberg fixed point theorem  \cite{ab} 
and $n$ Bellman's optimality equations for  dynamic programming problems  associated with the players. 
The literature on  dynamic programming in Markov decision models is very well described in  \cite{bs,bl0,hll,w}.

A number of  natural examples of static constrained games come from economics \cite{debreu,rosen}.  
Constrained $n$-person Markov games with finite state and action spaces were first studied by Altman and Shwartz \cite{as}.
They also apply in their analysis the Kakutani fixed point theorem and formulate the problem 
for each player as a Markov decision  process.
However, by introducing stochastic constraints on the strategy choices of the players, 
they had to use some facts from the theory 
of sensitivity analysis of  linear programming. 

Nash equilibria in games with constraints arise quite naturally, for instance, in the context of  
asynchronous transfer mode  networks, 
where users express their requirements for quality of service by bounds they wish to have on delays, etc. 
An audio application could therefore selfishly seek to minimise losses, 
subject to a maximum bound on the delay  it experiences.   
Nash equilibria in constrained games were also  studied   in a dynamic environment in telecommunications 
and internet provisioning applications, see \cite{1} and reference cited therein. 
Other applications focus on selection of rate allocations in multiple access channels
as well as models with asymmetric or partial information \cite{3,4}.  
Applications of constrained stochastic games to some queueing models are given in \cite{ahl,zhg}. 

The result of Altman and Shwartz \cite{as} for discounted constrained Markov games
was generalised by Alvarez-Mena and Hern{\'a}ndez-Lerma  \cite{ahl},
who considered compact metric action spaces. They considered first finite state space games and 
next, imposing a special condition on the transition probability, showed
how to get a stationary Nash equilibrium in a discounted game with countably many states. 
The proof of this result   in \cite{ahl}  relies on an approximation of the
game with denumerable state space by games with finitely many states. 

Work \cite{ahl} is devoted to games with bounded cost (payoff) functions.
The main result in \cite{ahl} was further used by Zhang et al. \cite{zhg} to prove the existence
of stationary Nash equilibria in a class of discounted Markov games with countable state spaces and unbounded cost functions. 
The assumptions, that Zhang et al. made in \cite{zhg}, resemble conditions presented by Wessels \cite{w}, 
who  studied dynamic programming 
problems with unbounded reward functions using the so-called weighted norm. Similarly as in \cite{ahl}, Zhang et al. \cite{zhg} apply an 
approximation of the original game by ones with finitely many states.

The value iteration algorithms and Bellman's optimality equations are
not sufficient tools for studying constrained Markov decision  problems and
constrained stochastic games.  As shown by Borkar \cite{bor}, some results from convex analysis and properties
of so-called  occupation measures induced by strategies (control functions)
must be applied. The approach using occupation measures enables to recognise the dynamic optimisation  problem 
as a static one on some compact convex subset of probability measures on the space of state-action pairs. 
The compactness and convexity of the set of occupation measures in the discounted Markov decision process   
is is closely related to the properties of the space of all probability measures on the set  of trajectories of the process 
induced by strategies (policies) of the decision maker.  For the details the reader is referred to \cite{bor,auto,piu}. 

The existence of stationary Nash equilibria in constrained discounted Markov games
is proved by using the Kakutani-Fan-Glicksberg fixed point theorem \cite{ab}. However, the main obstacle is 
to define the domain for the best-response correspondence associated with an auxiliary one-shot game.
In the {\it finite state space case}, studied by Altman and Shwartz \cite{as} and 
Alvarez-Mena and Hern{\'a}ndez-Lerma  \cite{ahl}, for any player $i$,
the authors take into account the set, say   $\Pr(\mathbb{K}_i),$ of all probability measures on the set  $\mathbb{K}_i$ 
of  all pairs $(x,a_i),$ where $x$ is a state and $a_i$ is an action available to player $i$  in this state. 
Then, $\Pr(\mathbb{K}_i)$ is convex and compact in the weak topology. The best responses
of   player $i$  also belong  to $\Pr(\mathbb{K}_i) $ and satisfy some equations introduced
in the theory of constrained Markov decision  models by Borkar \cite{bor}.
These equations guarantee that the fixed point of the best-response correspondence
is a vector of occupation measures, from which the existence of a stationary Nash equilibrium is concluded.  
When the state space is countable and {\it infinite}, then the spaces $\Pr(\mathbb{K}_i)$ 
are not compact and the analysis from \cite{as} and \cite{ahl} does not work. Therefore, 
Alvarez-Mena and Hern{\'a}ndez-Lerma  \cite{ahl} and Zhang et al. \cite{zhg}  
approximate  the stochastic game with denumerable state spaces by games with finitely many  states.  
Analogous methods were earlier used to consider discounted Markov decision  processes by Altman \cite{a} and 
Cavazos-Cadena \cite{cc}.  We would like to emphasise that introducing 
constraints in the stochastic game  model and following the finite state approximations as
in \cite{ahl,zhg} lead  to several unnecessary technical considerations. Therefore, our techniques and ideas 
are different than in the aforementioned papers.

In this paper, we study a general class of discounted  constrained 
stochastic  Markov games with unbounded costs. In Section 2, we   
formulate our basic assumptions including the uniform integrability of discounted  cost functions on the space of
all trajectories (sample paths) of the process. 
Our assumptions are weaker than those used by Zhang et al. \cite{zhg} for games and by Wessels \cite{w} for dynamic programming.  
Section 2  also presents our main results and  contains a few essential comments (Remarks 1-4) and remarks on earlier works \cite{as,ahl,zhg}. 
The proof of the main theorem (Theorem 1) for unbounded Markov games with stochastic constraints is provided in Section 4.
It is based on an auxiliary result 
(Proposition 1  in Section 3)  
for constrained Markov games with bounded costs and positive initial state distribution. An   
 approximation of the general stochastic  game by ones with  
 perturbed initial state distributions and truncated costs is applied.   	 
In Section 5, we give examples  that explain  relations of our uniform integrability  assumption 
from Section 2 with those of  Wessels \cite{w} and \cite{zhg}.
Section 6 (Appendix) contains two lemmas used in the proofs  in Sections 3 and 4.

Finally, we wish to stress out that our 
idea applied  for the study of  games 
with bounded costs is new and relies on introducing a proper domain for the best-response
correspondence associated with the auxiliary one-shot game. 
This domain is the Cartesian product of some 
appropriately constructed compact convex subsets of the spaces $\Pr(\mathbb{K}_i)$ 
(all probability measures on $\mathbb{K}_i$). 
Our approach works in the infinite countable state space case and therefore, no finite state approximation
is necessary. Instead, we apply the basic results on occupation 
and strategic measures from Borkar \cite{bor}
and Sch{\"a}l \cite{s1,s2}. 

\section{ The model and the main results}

The non-zero-sum {\it constrained stochastic Markov game} {\it (CSG)} is described by the following objects: 
\begin{itemize}
\item ${\cal N}=\{1,2,...,n\}$  is the {\it set of players}. 
\item $X$ is a countable {\it state space} endowed with the discrete topology. 
\item $A_i$ is a Borel {\it action space} for
player $i\in {\cal N}.$ The set $A_i(x)$ is a non-empty compact subset of $A_i,$ $x\in X,$ $i\in {\cal N}.$ We put
$$ A:=\prod_{i=1}^n A_i \quad\mbox{and}\quad 
 A(x):= \prod_{i=1}^n A_i(x).$$
Note that the set
$$\mathbb{K}_i= \{(x,a_i): x\in X,\ a_i\in A_i(x) \}$$
of {\it feasible state-action pairs} for player $i\in{\cal N}$ is a closed subset of $X\times A_i$. Similarly, the set 
$$\mathbb{K}= \{(x,\pmb{a}): x\in X,\ \pmb{ a}=(a_1,...,a_n)\in A(x) \}$$
of {\it feasible state-action vectors} is a closed subset of $X\times A.$
\item Let $L=\{1,...,l\}$ and $L_0=L\cup\{0\}.$ The real-valued functions
$c_i^\ell:\mathbb{K}\to \mathbb{R},$ $i\in {\cal N},$ $\ell \in L_0$ are  measurable. Here,
$c_i^0$ denotes {\it cost-per-stage  function} for player  $i\in {\cal N},$ 
and for each  $\ell\in L,$   $c_i^\ell$ is a function used in the definition of the $\ell$-th {\it constraint} for this player. 
\item $p(y|x,\pmb{a})$ is the transition probability from $x$ to $y\in X,$
when the players choose a profile $\pmb{ a}=(a_1,a_2,...,a_n)$ of actions in $ A(x).$ 
\item $\eta$ is  the {\it initial state distribution}.
\item $\alpha\in (0,1)$ is the {\it discount factor}.
\item $\kappa_i^\ell$ are constants, $i\in {\cal N},$ $\ell \in L.$ 
\end{itemize}
Let $\mathbb{N}=\{1,2,...\} .$ Define $H^1=X$ and $H^{t+1}=   \mathbb{K}\times H^{t}$  for $t\in \mathbb{N}.$ 
An element  $h^t=(x^1,\pmb{ a}^1,\ldots,x^t)$ of $H^t$ represents a history of the game up to the $t$-th
stage, where  $\pmb{ a}^k=(a^k_1,\ldots,a^k_n)$ is the profile of actions chosen by the players in the state $x^k$ 
on the $k$-th stage of the game 
($k\in \mathbb{N}$). Clearly, $h^1=x^1.$

Strategies for the players are defined in he usual manner. 
A {\it strategy} for player $i\in{\cal N}$ is a sequence $\pi_i=(\pi_{i}^t)_{t\in  \mathbb{N}},$ 
where each $\pi_{i}^t$  is a transition probability from $H^t$ to $A_i$ such that
$\pi_{i}^t(A_i(x^t)|h^t)=1$ for any history $h^t\in H^t,$ $t\in  \mathbb{N}.$
By $\Pi_i$ we denote the {\it set of all strategies} for player $i.$
Let  $\Phi_i$ be the set of transition probabilities from $X$ to $A_i.$ Then,
$\varphi_i\in\Phi_i$ if $\varphi_i(A_i(x)|x)=1$ for all $x\in X.$ A {\it stationary strategy} for player $i$ is a constant
sequence $ (\varphi_{i}^t)_{t\in  \mathbb{N}},$ where $\varphi_i^t=\varphi_i$ for all $t\in\mathbb{N}$ and some $\varphi_i\in\Phi_i.$ 
Furthermore, we shall identify  a stationary strategy for player $i$  with the  constant element $\varphi_i$  of the sequence. 
Thus, the {\it set of all    stationary } strategies of player $i$  is $\Phi_i.$
We define
$$\Pi = \prod_{i=1}^n \Pi_i \quad\mbox{and}\quad
\Phi = \prod_{i=1}^n \Phi_i.$$
Hence, $\Pi$ ($\Phi$) is the set of all (stationary) multi-strategies of the players.

Let  $H^\infty = \mathbb{K}\times \mathbb{K}\times \cdots$
be the space of all infinite histories of the game endowed with
the  product $\sigma$-algebra. For any  multi-strategy  $\pmb{ \pi}\in\Pi$,
a probability measure $\mathbb{P}_\eta^{\pmb{\pi}}$  and a stochastic process $(x^t,\pmb{a}^t)_{t\in \mathbb{N}}$  
are defined on $H^\infty$ in a canonical way, see the Ionescu-Tulcea theorem, e.g., Proposition 7.28 in \cite{bs}. 
The measure $\mathbb{P}_\eta^{\pmb{\pi}}$  is induced by $\pmb{\pi},$ the transition probability $p$ and the initial distribution $\eta.$  
The expectation operator with respect to $\mathbb{P}_\eta^{\pmb{\pi}}$ is denoted by $\mathbb{E}_\eta^{\pmb{\pi}}.$ 

Let $\pmb{\pi}\in\Pi$ be any multi-strategy. For each
$i\in {\cal N}$ and $\ell \in L_0$,  the {\it discounted cost functionals} (see   \cite{bs,bl0})   are     defined as follows:
\begin{equation}\label{jil}
J_i^\ell(\pmb{\pi}) = (1-\alpha)\mathbb{E}^{\pmb{\pi}}_\eta\left[\sum_{t=1}^\infty  \alpha^{t-1} c_i^\ell(x^t,\pmb{a}^t) \right].
\end{equation}
Below we provide conditions that guarantee that the functionals are well-defined. 
 We assume that $J^0_i(\pmb{\pi})$ is the expected discounted cost of player $i\in {\cal N}$,  who wishes to   
minimise it over $\pi_i \in\Pi_i$ in such a way  that
the following constraints are satisfied
\begin{equation}
\label{ico}
J^\ell_i(\pmb{\pi}) \le \kappa^\ell_i \quad\mbox{for all}\quad
\ell\in L.
\end{equation}
A multi-strategy $\pmb{\pi}$ is {\it feasible}, if (\ref{ico})
holds for each $i\in {\cal N},$  $\ell\in L.$ We denote by $\Delta$ the set of all feasible multi-strategies in {\it CSG.}

As usual, for any $\pmb{\pi}\in\Pi$,  we denote by  $\pmb{\pi_{-i}}$
the multi-strategy of all players but player $i.$  More precisely,
$ \pmb{\pi_{-1}} =(\pi_2,...,\pi_n),$  $ \pmb{\pi_{-n}} =(\pi_1,...,\pi_{n-1}),$
and for $i\in {\cal N}\setminus \{1,n\},$  
 $$\pmb{\pi_{-i}}
=(\pi_1,\ldots,\pi_{i-1},\pi_{i+1},\ldots,\pi_n).$$
We identify $(\pmb{\pi_{-i}},\pi_i)$ with $\pmb{\pi}.$ 
 For each $\pmb{\pi}\in\Pi$, we define the set of feasible strategies for player $i$ with $\pmb{\pi_{-i}}$ as
$$\Delta_i(\pmb{\pi_{-i}})=\{\pi_i\in\Pi_i: \ J^\ell_i(\pmb{\pi}) 
=  J^\ell_i(\pmb{\pi_{-i}},\pi_i) 
\le \kappa^\ell_i \quad\mbox{for all}\quad
\ell\in L\}.$$
Hence, $\pmb{\pi}\in\Delta$ if and only if $\pi_i\in\Delta_i(\pmb{\pi_{-i}})$ for all $i\in {\cal N}.$

Let $\pmb{\pi}=(\pi_1,\pi_2,...,\pi_n) \in\Pi$ and $\sigma_i\in\Pi_i.$   By $[\pmb{\pi_{-i}},\sigma_i]$ we denote the multi-strategy, 
where player $i$ uses $\sigma_i$ and every player $j\not= i$ uses $\pi_j.$\\

\noindent{\bf Definition 1.} A multi-strategy $\pmb{\pi}\in\Pi$ is a {\it Nash equilibrium} in {\it CSG}, if 
$\pi_i\in \Delta_i(\pmb{\pi_{-i}})$    and 
$$J^0_i(\pmb{\pi})=\inf_{\sigma_i \in\Delta_i(\pmb{\pi_{-i}})} J^0_i([\pmb{\pi_{-i}},\sigma_i])$$ for every player $i\in {\cal N}.$
 \\

\noindent{\bf  Assumption A} \\{\it 
{\bf (i)}  The function $p(y|x,\cdot)$ is continuous on $A(x)$ for all $x,\ y\in X.$\\
{\bf (ii)} The functions $c^\ell_i(x,\cdot)$ are continuous on $A(x)$ for all $x\in X,$ $i\in {\cal N}$ and $\ell\in L_0.$}
\\

\noindent{\bf  Assumption B} \\{\it 
{\bf (i)} There exists a  function $w: X\to[1,\infty)$
such that $|c_i^\ell(x,\pmb{a})|\le w(x)$ for each $(x,\pmb{a})\in\mathbb{K} $ and  for all $\ell\in L_0,$ $i\in{\cal N}.$ \\
{\bf (ii)} It holds
\begin{equation}\label{Bii1}
 \lim_{k\to\infty} \sup_{\pmb{\pi}\in\Pi}  (1-\alpha) \mathbb{E}_\eta^{\pmb{\pi}}\left[\sum_{t=k}^\infty \alpha^{t-1}w(x^t) \right]=0,
\end{equation}
and, for each $t\in\mathbb{N},$
\begin{equation}\label{Bii2}
 \lim_{k\to\infty} \sup_{\pmb{\pi}\in\Pi} \mathbb{E}_\eta^{\pmb{\pi}}\left[w(x^t)1_{[w(x^t)\ge k]} \right]=0,
\end{equation}
where $1_K$ denotes the indicator function of the set $K$. 
}
\\

A weaker version of assumption {\bf B(i)} was used in \cite{auto} to study constrained Markov 
decision processes on a Borel state space. 
Assumption {\bf B} implies that all expectations in (\ref{jil})  are finite. To show this fact 
fix $\epsilon >0.$ From (\ref{Bii1}),
there exists $N_0\in\mathbb{N}$  such that for all $k>N_0,$ we have
\begin{eqnarray}\label{ogr}\nonumber
|J_i^\ell(\pmb{\pi}) |&\le&  (1-\alpha)\mathbb{E}^{\pmb{\pi}}_\eta\left[\sum_{t=1}^\infty  \alpha^{t-1} w(x^t) \right] \nonumber
=  (1-\alpha)\mathbb{E}^{\pmb{\pi}}_\eta\left[\sum_{t=1}^{k}  \alpha^{t-1} w(x^t) \right]\\ &+& 
(1-\alpha)\mathbb{E}^{\pmb{\pi}}_\eta\left[\sum_{t=k+1}^\infty  \alpha^{t-1} w(x^t) \right] 
 \le  (1-\alpha)\mathbb{E}^{\pmb{\pi}}_\eta\left[\sum_{t=1}^{k}  \alpha^{t-1} w(x^t) \right] +  \epsilon.
\end{eqnarray}
for all $\pmb{\pi}\in \Pi$ and $i\in{\cal N}$ and $\ell\in L_0.$ Now let us consider the first term on the right-hand side in (\ref{ogr}). 
By (\ref{Bii2}), we may choose $m\in\mathbb{N}$
such that  for each $t=1,\ldots, k$ we obtain the bound
\begin{eqnarray*}
 (1-\alpha)\mathbb{E}^{\pmb{\pi}}_\eta\left[\sum_{t=1}^{k}  \alpha^{t-1} w(x^t) \right] 
&\le&  
(1-\alpha)\left[\sum_{t=1}^{k}  \alpha^{t-1} m\right] +(1-\alpha)\mathbb{E}^{\pmb{\pi}}_\eta\left[\sum_{t=1}^{k}  \alpha^{t-1} w(x^t) 1_{[w(x^t)\ge m]} \right]\\
&<&m+ (1-\alpha)\left[\sum_{t=1}^{k}  \alpha^{t-1} \epsilon\right]
<m+\epsilon.
\end{eqnarray*}
Summing up, we get that 
$$|J_i^\ell(\pmb{\pi}) |< m+2\epsilon<\infty.$$
Similar  assumptions to study  non-stationary Markov decision processes with unbounded payoffs were  formulated and thoroughly discussed
in \cite{dg}.
Further details and comments,   the reader   can find in Section 5.

The next assumption  is called
in the literature the {\it Slater condition}, see \cite{as},
Assumption 3.3(c) in \cite{ahl} and Assumption 2 in \cite{zhg}. \\

\noindent{\bf  Assumption C} \\ {\it For each stationary multi-strategy
$\pmb{\varphi}\in \Phi$ and for each player $i\in {\cal N},$ there exists
$\pi_i\in \Pi_i$ such that}
$$
J^\ell_i([\pmb{\varphi_{-i}},\pi_i]) < \kappa^\ell_i  \ \mbox{for all}\ \ell \in L.
$$

Observe that under assumption  {\bf C} the set of feasible strategies  in $CSG$ is  non-empty. \\

We are ready to state our main result. \\

\noindent{\bf Theorem 1.} {\it Assume {\bf A}, {\bf B} and {\bf C}. Then, the CSG possesses a Nash equilibrium in the set $\Phi.$}\\

The proof of this result is given in Section 4.

Below we describe a special case of our Assumption {\bf  B}. \\

\noindent{\bf  Assumption W} \\{\it 
{\bf (i)} There exists a  function $w: X\to[1,\infty)$
such that {\bf B(i)} holds and 
$$\sum_{y\in X} w(y)p(y|x,\pmb{a})\le \delta w(x)\quad\mbox{for all}\ \  
(x,\pmb{a})\in\mathbb{K},$$
and for some constant $\delta\ge 1$ satisfying $\delta\alpha<1.$\\
{\bf (ii)} The function $\sum_{y\in X} w(y)p(y|x,\cdot)$ is continuous on $A(x)$ for each $x\in X.$\\
{\bf (iii)} $\sum_{y\in X} w(y)\eta(y) < \infty.$ }
 \\

The above conditions were first  introduced in \cite{w} to deal with unbounded payoffs in Markov decision processes. 
They gained recognition and 
were broadly applied to deal with several models, see  for instance \cite{hll,jn}.

By  Lemma 9  in \cite{auto}, it follows that, if  there exists a function $w$ that satisfies {\bf W}, then {\bf B} holds as well.
However, Example 4 in \cite{auto} warns that the violation of {\bf W(ii)} entails that the uniform integrability condition in
(\ref{Bii2}) fails. Moreover, this example (case II on p. 10 in \cite{auto}) also illustrates 
that,  if {\bf B} is satisfied with any value of a discount factor, then 
{\bf W(i)} holds only for $\alpha<4/5.$

From Theorem 1 we can deduce two conclusions.\\

\noindent{\bf Corollary 1.} {\it Assume {\bf A}, {\bf W} and {\bf C}. Then, the CSG possesses a Nash equilibrium in the set $\Phi.$}\\

\noindent{\bf Corollary 2.} {\it Assume {\bf A},   {\bf C}
and that every function $c_i^\ell$ is bounded for $i\in {\cal N},$ $\ell\in L_0.$
Then, the CSG possesses a Nash equilibrium in the set $\Phi.$}\\

\noindent{\bf Remark 1.} 
Discounted constrained stochastic games with countable state spaces
and unbounded functions $c^\ell_i$ were studied by Zhang et al. \cite{zhg}.  
However, the assumptions imposed  in \cite{zhg} are stronger than ours. Indeed, they require that there exists an unbounded function $w$ that
satisfies   {\bf B(i)} and such that $w^2$ is integrable with respect to the initial state  distribution $\eta$  and with respect to the transition probability. 
More precisely, Assumption 1(e) in \cite{zhg}  says that
there exists a constant $\beta\ge 1$ such that $\alpha\beta^2<1$ and 
\begin{equation}
\label{zhang}
\sum_{y\in X}  w^2(y)p(y|x,\pmb{a})\le \beta^2  w^2(x)\quad
\mbox{ for all} \quad (x,\pmb{a})\in\mathbb{K}.
\end{equation}
In addition, $w$ is  a moment function, i.e., there exists an increasing sequence of finite sets $(Z_m)_{m\in\mathbb{N}}$ such that
$\bigcup_{m\in\mathbb{N}} Z_m =X$ 
and  $\lim_{m\to\infty} \inf_{x\not\in Z_m} w(x)=\infty.$   This condition excludes from consideration  games on infinite countable 
state spaces with bounded cost functions and arbitrary transition probabilities, since condition (\ref{zhang}) 
need not hold for the unbounded function $w.$  

Note that   (\ref{zhang}) can be written in the form
$$
\sum_{y\in X}  \widehat{w}(y)p(y|x,\pmb{a})\le \delta  \widehat{w}(x)\quad
\mbox{ for all} \quad (x,\pmb{a})\in\mathbb{K} 
$$
where $\widehat{w}= w^2$, $\delta = \beta^2, $ $\delta\alpha <1.$ This is exactly condition
{\bf W(i).} But then {\bf B(i)} (assumed in {\bf W(i)}) holds as well, that is, 
\begin{equation}
\label{linout}
|c_i^\ell(x,\pmb{a})|\le \sqrt{\hat{w}(x)}\quad\mbox{ for all} \quad
(x,\pmb{a})\in\mathbb{K},\    \ell\in L_0, \ i\in{\cal N}.
\end{equation} 
Thus, the class of games satisfying our assumption {\bf W}
is essentially larger than the class studied in \cite{zhg}.
For example, condition (\ref{linout}) excludes linear functions $c_i^\ell,$
when $\widehat{w}$ is linear
as it happens in many Markov decision processes, for examples consult  with \cite{auto,jn}.  
Therefore, Corollary 1 extends Theorem 1 in \cite{zhg}. \\

\noindent{\bf Remark 2.} Discounted 
constrained stochastic games with finite state and action spaces were first studied by Altman and Shwartz   \cite{as}. 
An extension to games with compact metric action spaces  was  given by Alvarez-Mena and 
Hern{\'a}ndez-Lerma  \cite{ahl}. The existence of stationary Nash equilibrium 
in the finite state space framework with constraints is established by a fixed point argument, 
but the approach from the unconstrained case as in \cite{f,fink,taka} cannot be applied. 
The main difficulty is to determine a domain for the best-response correspondence, sometimes  
called the Nash correspondence. Unlike the standard case \cite{f,fink,taka}, the Cartesian product of the sets 
of stationary strategies is not appropriate in the constrained setting. The authors consider the Cartesian product of the sets $\Pr(\mathbb{K}_i)$ 
of all probability measures on $\mathbb{K}_i$ ($i\in {\cal N})$ and an auxiliary one-shot game. The sets $\Pr(\mathbb{K}_i)$ 
are actually too large and, therefore, some functional equations, characterising so-called ``occupation measures'' on $\mathbb{K}_i,$  are 
requested in the definition of the Nash correspondence.
These equations 
play a fundamental role in discounted constrained decision processes
and games \cite{ahl,bor,piu}.  In the finite state space case, the sets $\Pr(\mathbb{K}_i)$ 
are compact in the weak topology and obviously they are convex. 
Consequently, the Kakutani-Fan-Glicksberg fixed point theorem \cite{ab} can be applied. 
When the state space $X$ is countable and infinite, then
the Cartesian product of the spaces $\Pr(\mathbb{K}_i)$ cannot be used as a domain 
for the Nash correspondence, because $\Pr(\mathbb{K}_i)$ is
non-compact in the weak topology. Therefore, to study discounted $CSGs$ 
with an infinite countable state space $X$,     
Alvarez-Mena and  Hern{\'a}ndez-Lerma  \cite{ahl} and Zhang et al. \cite{zhg} use 
an approximation of the original game by games with finite state spaces. 
The proof in \cite{ahl}  strongly exploits  Assumption 3.4 (see p. 267 in \cite{ahl}) that entails 
the convergence of discounted costs in the approximating models to the discounted cost in the original model  
(consult with the proof  Theorem 3.6(c) in \cite{ahl}).
Using our notation,  Assumption 3.4 from \cite{ahl} sounds as follows: there exists 
an increasing sequence of finite sets $(Z_m)_{m\in\mathbb{N}}$ such that
$\bigcup_{m\in\mathbb{N}} Z_m =X$  and  
\begin{equation}
\label{ahlass}
\lim_{m\to\infty} \max_{x\in Z_m}\max_{\pmb{a}\in A(x)}
p(X\setminus Z_{m+1}|x,\pmb{a}) =0.
\end{equation}
The cost and constraint functions in \cite{ahl}, however, are bounded
and additionally, the condition in (\ref{ahlass}) looks restrictive. 
Zhang et al. \cite{zhg} also 
approximate the original discounted $CSG$  by appropriately defined  auxiliary finite state space games 
and show that their stationary Nash equilibria   converge to a   Nash equilibrium    
in the original game. As mentioned in Remark 1, they allow the functions $c^\ell_i$ to be unbounded and 
drop Assumption 3.4 from \cite{ahl}. Their proof, on the other hand, is inspired  
by an estimation techniques developed in \cite{cc} and Chapter 16 in \cite{a}.\\ 

\noindent{\bf Remark 3.} 
The proof of Theorem 1 proceeds along   different lines than those in \cite{ahl,zhg}. 
First of all, we do not apply the finite state space approximations. Our proof is more direct. 
The idea is based  on studying first auxiliary stochastic games  with bounded cost functions and with  
positive initial state distributions.  
An important new feature of our approach is to define the Nash correspondence using  
some compact convex subsets of the spaces $\Pr(\mathbb{K}_i)$ 
being projections of the  set $Y$ of occupation measures induced by correlated strategies of the players.  
The set $Y$ is a compact convex subset of the non-compact set of all probability measures on $\mathbb{K}.$ 
This idea  combined with basic results from convex analysis in Markov decision  processes, see., e,g.,  \cite{bor,piu}, 
establishes in Section 3 the existence of stationary Nash equilibria in  the   auxiliary discounted $CSGs.$
In Section 4, we show how to approximate the original game by the aforementioned auxiliary games with bounded costs. 
We prove that there is a sequence of Nash  equilibria 
in auxiliary games converging to a Nash equilibrium in the original discounted $CSG.$ 
More detailed comments on the basic idea in this paper used in studying $CSGs$ with bounded costs  
are given in Remark 5 in Section 3.  \\ 

\noindent{\bf Remark 4.} 
Corollary 2 states that there exists a stationary Nash equilibrium in games with bounded cost and constraint  functions.  
The same result is formulated as Corollary 1 in \cite{zhg}. However,  the proof of Corollary 1 in \cite{zhg} is incorrect. 
Firstly, it cannot be deduced from Theorem 1 in their paper, where $w$ is assumed to be a  moment function. 
Using an unbounded function $w$ one has to restrict the classes of transition
probabilities and initial distributions (e.g., the integral of $w$ with respect to $\eta$ need not be finite). 
Secondly,  in the proof of Corollary 1 in \cite{zhg}, Zhang et al. 
erroneously claim that every transition probability satisfies
condition (\ref{ahlass}) (made as Assumption 3.4 in \cite{ahl})
with $Z_m=\{1,2,...,m\}.$   That is not true.  For example, (\ref{ahlass})
fails to hold when  $p(x+2|x,\pmb{a})=1$ for all $x\in X$ and $\pmb{a}\in A(x).$   
Claiming that (\ref{ahlass}) holds in general, Zhang et al. \cite{zhg}  
conclude their Corollary 1 from Theorem 3.6(c) in \cite{ahl}.

 \section{Stochastic games  with bounded costs and positive state distribution}

In this section, we state an auxiliary result  (Propostion 1)  using basic theorems on occupation and strategic measures 
obtained by Borkar \cite{bor} and  Sch{\" a}l \cite{s1,s2}. 

\subsection{Occupation measures and their important properties}

Let $Y$ be a metric space with the Borel $\sigma$-algebra ${\cal B}(Y).$ Let $\Pr(Y)$ ($\mathbb{M}(Y)$) be the set 
of all probability (finite signed) measures on $Y$
and $C(Y)$ be the space of  all bounded uniformly continuous functions on $Y$. A sequence $(\nu^k)_{k\in\mathbb{N}}$  in $\mathbb{M}(Y)$ is said to
{\it converge weakly} to $\nu\in\mathbb{M}(Y)$ if 
$$\int_Y fd\nu^k\ \to\ \int_Yfd\nu\quad \mbox{for all} \ f\in C(Y).$$
If $Y$ is compact metric, then so is $\Pr(Y)$ equipped with the topology of weak convergence, see  Theorem 6.4 in \cite{par}. 
If $Y$ is a Borel space, then we may equip $\mathbb{M}(Y)$ with the metric
$$d(\mu,\nu)=\sum_{k=1}^{\infty}\frac 1{2^k}\left|\int_Y f_kd\mu-\int_Y f_kd\nu \right|, \quad \mu,\ \nu\in\mathbb{M}(Y).$$
Here, ${\cal F}=\{f_k\}_{k\in\mathbb{N}}$ is a countable family of
functions, which is dense in the unit ball in $C(Y)$ and such that for each different points $y,\ y'\in Y,$
there exists $f_k\in{\cal F}$ such that $f_k(y)\not= f_k(y')$ (see p. 47 in \cite{par}).
It is obvious that the topology induced by this metric is equivalent to the weak topology in 
$\mathbb{M}(Y)$. Thus, $\mathbb{M}(Y)$ is a Hausdorff locally convex space.

Since $A_i(x)$ is compact metric for each $x\in X,$   so is $\Pr(A_i(x))$ 
in the weak topology. 
The set $\Phi_i$ can be identified with the product space
$\prod_{x\in X} \Pr(A_i(x)).$
By the Tychonoff theorem the spaces
$$\Phi_i=\prod_{x\in X} \Pr(A_i(x))\quad\mbox{and}\quad \Phi=\prod_{i\in {\cal N}} \Phi_i$$ 
are compact, when endowed with the product topologies. 
Moreover, these spaces are metrisable.

A sequence $(\varphi_i^{k})_{k\in\mathbb{N}}$ in $\Phi_i$ converges   to $\varphi_i\in\Phi_i,$ if the sequence 
$(\varphi_i^k(\cdot|x))_{k\in\mathbb{N}}$ in $\Pr(A_i(x))$ converges weakly to $\varphi_i(\cdot|x)$ for each $x\in X.$
A sequence $(\pmb{\varphi}^k)_{k\in\mathbb{N}}$ in $\Phi$ converges   to $\pmb{\varphi}\in \Phi,$ if 
$(\varphi_{i}^k)_{k\in\mathbb{N}}$ converges   to $\varphi_i$ in $\Phi_i$ for every $i\in{\cal N}.$

For some technical reasons we also introduce $\widetilde{\Pi}$
as the {\it set of all  correlated strategies} $\pi=(\pi^t)_{t\in \mathbb{N}}$ of the players.
Here, $\pi^t$ is a transition probability from $H^t$ to $A$ such that
$\pi^t(A(x^t)|h^t)=1$ for any history $h^t\in H^t.$
Using a correlated strategy the players act like  one decision maker in the Markov decision process
with the action spaces $A(x)$, $x\in X.$

Let $\cal M$ be the set of probability {\it occupation measures} on $\mathbb{K}$
induced by all correlated strategies $\pi \in\widetilde{\Pi}$
and the initial distribution $\eta,$ i.e., $\rho\in {\cal M}$ is defined as follows 
\begin{equation}
\label{ro0}
\rho(K)= (1-\alpha)\mathbb{E}_\eta^\pi\left[\sum_{t=1}^\infty \alpha^{t-1} 1_K(x^t,\pmb{a}^t)\right]\quad\mbox{
for  }K\in{\cal B}(\mathbb{K}).
\end{equation}
The expectation operator  $\mathbb{E}_\eta^\pi$ is taken with respect to the unique probability measure
 $\mathbb{P}_\eta^{\pi}$ on $H^\infty$, called a {\it strategic measure}.
 
The integral of any bounded measurable function $f:A(x)\to \mathbb{R}$
with respect to $\rho(\{x\}\times \cdot)$ is denoted by
$\int_{A(x)}f(\pmb{a})\rho(\{x\}\times d\pmb{a}).$

From (\ref{ro0}), it follows that for any bounded measurable function $c:\mathbb{K}\to \mathbb{R}$,
\begin{equation}
\label{ro1}
\int_\mathbb{K}cd\rho = \sum_{x\in X} \int_{A(x)}c(x,\pmb{a})\rho(\{x\}\times d\pmb{a}) =
(1-\alpha)\mathbb{E}_\eta^\pi\left[\sum_{t=1}^\infty \alpha^{t-1} c(x^t,\pmb{a}^t)\right].
\end{equation}
Hence, for any $i\in {\cal N}$,  $ \ell \in L_0 $  and $\pi\in \widetilde{\Pi},$ 
$$\int_\mathbb{K}c_i^\ell d\rho = 
\sum_{x\in X} \int_{A(x)}c^\ell_i(x,\pmb{a})
\rho(\{x\}\times d\pmb{a})
= J_i^\ell(\pi).
$$
\\
 
\noindent{\bf  Lemma 1. } {\it Under assumption {\bf A(i)}}, $\cal M$ {\it is  convex 
and   compact  in $\Pr(\mathbb{K})$ equipped with the weak topology.} \\

{\it Proof.}  The set $\{ \mathbb{P}_\eta^{\pi}: \pi \in 
\widetilde{\Pi}\}$ of all strategic measures induced by all correlated strategies is weakly compact in $\Pr(H^\infty).$  
It is also convex.
These facts  are well-known in the literature and 
together with (\ref{ro1}) imply the lemma. For a detailed discussion see: \cite{piu}, Subsections 7.1 and 7.4 in \cite{s2}, 
Theorem 5.6 in \cite{s1}, or  Theorem 3.1 in \cite{bor}.  	\hfill $\Box$\\

Let  ${\cal M}_i$ be the set of all probability measures on $\mathbb{K}_i$ defined as follows. 
A measure $\mu$ belongs to ${\cal M}_i$, if there exists a probability measure $\rho \in {\cal M}$ 
such that, for each $x\in X,$  $\mu(\{x\}\times \cdot)$ is the projection of $\rho(\{x\}\times \cdot)$ on $A_i(x).$
More detailed, if $x\in X,$   $B\in {\cal B}(A_i(x))$
and $pr$ is the projection from $A(x)$ on $A_i(x)$, then
$$
\mu(\{x\}\times B)= \rho(\{x\}\times  pr^{-1}(B)).
$$
If $f$ is a bounded measurable real-valued function on $A_i(x),$
then $\int_{A_i(x)} f(a_i) \mu(\{x\}\times da_i)$   means the integral of $f$
with respect to the probability measure 
$ \mu(\{x\}\times \cdot).$ 
 Every   function $f\in C(\mathbb{K}_i)$ 
 can be recognised as a  function in $C(\mathbb{K}).$ Then
\begin{equation}
\label{pr}
\int_{\mathbb{K}} f d\rho=  \int_{\mathbb{K}_i} f d\mu
= \sum_{x\in X}\int_{A_i(x)}f(x,a_i)\mu(\{x\}\times da_i).
\end{equation}
Hence, if $(\rho^k)_{k\in \mathbb{N}}$  is a sequence of measures
converging weakly to some $\rho^0$    in $\cal M$  and $\mu^k,$ $\mu^0$ are
projections of $\rho^k$ and $\rho^0$, respectively, defined as above, then
by (\ref{ro1}) and  (\ref{pr}), 
$(\mu^k)$  converges weakly to $\mu^0$ in ${\cal M}_i.$
This fact and convexity of the set $\cal M$ imply the following result.\\

\noindent{\bf  Lemma 2.}   {\it If {\bf A(i)} holds, then}  ${\cal M}_i$ {\it is a convex and   compact  subset in} 
$\Pr(\mathbb{K}_i) $ {\it equipped with the weak topology.}\\

\noindent{\bf Remark 5.}  Introducing the sets $\widetilde{\Pi}$ and $\cal M$ is {\it crucial} in our proof. 
It enables us to use the compact and convex sets ${\cal M}_i$ in our definition of the best-response correspondence $S$ given below.  
In contrast to the finite state space case \cite{as,ahl},  $\Pr(\mathbb{K}_i)$ in games with  the infinite countable state spaces
cannot be utilised,  they need not be  compact. \\

In this subsection, we add the following condition.\\

\noindent{\bf Assumption D} \\
{\it For all $x\in X,$ $\eta(x)>0.$}\\

Let $\widehat{\mu}$ denote the projection of $\mu\in {\cal M}_i$ on $X,$
i.e., $\widehat{\mu}(x)= \mu(\{x\}\times A_i(x)), $  $x\in X.$\\

\noindent{\bf\ Lemma 3} {\it Assume {\bf A(i)} and  {\bf D}. If $\mu\in {\cal M}_i$, then 
$\widehat{\mu}(x)>0$   for all  $x\in X$    and there exists a unique  $\varphi_i\in\Phi_i$  such that 
\begin{equation}\label{dis}
\mu(\{x\}\times B)=\varphi_i(B|x)\widehat{\mu}(x),\quad \mbox{for all}
\quad    B\in{\cal B}(A_i(x)),\  x\in X. \end{equation}  }

{\it Proof.} 
Let $\mu\in{\cal M}_i$ be a projection of $\rho\in{\cal M}$
induced by some   $\pi\in \widetilde{\Pi}$ according to (\ref{ro0}). 
Then, for any $x\in X,$ 
\begin{eqnarray*}
\widehat{\mu}(x)&=&\mu(\{x\}\times A_i(x))=\rho(\{x\}\times A(x))= 
(1-\alpha)\sum_{t=1}^\infty \alpha^{t-1} \mathbb{E}_\eta^\pi\left(1_{\{x\}\times A(x)}(x^t,\pmb{a}^t)\right)\\
&\ge& 
(1-\alpha)\mathbb{P}_\eta^\pi(x^1=x)=(1-\alpha)\eta(x)>0.
\end{eqnarray*}
Therefore, $\varphi_i$ defined by 
$$ 
\varphi_i(B|x):= \frac{\mu(\{x\}\times B)}{\widehat{\mu}(x)}
\quad\mbox{for all}\quad  B\in{\cal B}(A_i(x)), \  x\in X,$$
is the unique transition probability satisfying  (\ref{dis}).
\hfill $\Box$
\\

\noindent{\bf  Remark 6.} Lemma 3 and assumption {\bf D} allow to omit the study of so-called equivalence classes 
of functions in $\Phi_i,$
which are equal on the set $Z\subset X$ with $\eta(Z)=1.$  They were
considered in \cite{ahl}. In our case,  
condition {\bf D} implies  the uniqueness of $\varphi_i\in\Phi_i$ in the above lemma and 
this fact simplifies our proofs in   the sequel.\\
 
Further,   we shall write  
(\ref{dis})  in the abbreviated form
$$\mu= \widehat{\mu}\varphi_i.$$
\\

\noindent{\bf  Lemma 4.} {\it 
Let $\mu^k= \widehat{\mu}^k\varphi_i^k,$ $k\in \mathbb{N}\cup
\{0\}.$ Under assumption {\bf A(i)} and {\bf D}, if $\mu^k \to \mu^0$ weakly in
${\cal M}_i,$ then,  for each $x\in X,$ 
$\widehat{\mu}^k(x) \to \widehat{\mu}^0(x)$ and
$\varphi_i^k(\cdot|x)\to \varphi_i^0(\cdot|x)$ weakly
in $\Pr(A_i(x))$ as $k\to\infty.$}\\

{\it Proof.}   
Since $\mu^k \to \mu^0$ weakly in
${\cal M}_i$, it follows that  $\widehat{\mu}^k(x)\to \widehat{\mu}^0(x)$ for every $x\in X.$  
Therefore, by Lemma 3, for every $x\in X,$   
$$\varphi_i^k(\cdot|x) = \frac{\mu^k(\{x\}\times\cdot)}{\widehat{\mu}^k(x)}\ 
\rightarrow\ \varphi^0_i(\cdot|x)= \frac{\mu^0(\{x\}\times\cdot)}{\widehat{\mu}^0(x)} \quad
 \mbox{in }\Pr(A_i(x))$$
endowed with the weak topology. 
\hfill $\Box$\\

\subsection{The existence of Nash equilibria in games with bounded costs}

In this subsection, we add the following assumption.\\

\noindent{\bf  Assumption B'}\\ {\it The functions 
$c_i^\ell$ are  bounded for all $i\in {\cal N}$ and $\ell\in L_0$.}\\

Under {\bf  B'} all functionals are bounded: $|J_i^\ell(\pmb{\pi})|\le \hat{c}$  
for some  $\hat{c}>0$ and for all $i\in {\cal N},$  $\ell \in L_0$  
and $\pmb{\pi}\in\Pi.$ 
Moreover, {\bf B(i)} and  {\bf B(ii)} are trivially satisfied
by taking $w(x)= \hat{c}$  for all $x\in X.$ 

Let $\pmb{\varphi}= (\varphi_1,\varphi_2,...,\varphi_n)\in \Phi.$ 
We denote by $[\pmb{\varphi_{-i}}(x),a_i]$ the profile of actions used in 
state $x\in X$, where player $i$ chooses $a_i\in A_i(x)$ and every player $j \not= i$ uses
$\varphi_j(\cdot|x).$   Then,
for $i\in {\cal N},$ $\ell \in L_0,$ $x\in X, $
$$
c^\ell_i(x, [\pmb{\varphi_{-i}}(x),a_i]):= \int_{A(x)}c_i^\ell(x,a_1,...,a_n)
\psi_1(da_1|x)...\psi_{n}(da_{n}|x)
$$
  and 
$$ p(y|x,  [\pmb{\varphi_{-i}}(x),a_i]):=
\int_{A(x)}p(y| x,a_1,...,a_n)
\psi_1(da_1|x)...\psi_{n}(da_{n}|x),$$
where $\psi_i(\{a_i\}|x)=1$ and $\psi_j =\varphi_j$ for all $j\not=i.$

To illustrate this notation, take $n=4,$ $i=3$ and 
$\psi_3(\{a_3\}|x)=1.$  Then, we get
$$
c^\ell_3(x, [\pmb{\varphi_{-3}}(x),a_3]) = \int_{A_4(x)}\int_{A_2(x)}\int_{A_1(x)}c_3^\ell(x,a_1,a_2,a_3,a_4)
\varphi_1(da_1|x)\varphi_2(da_2|x)\varphi_4(da_4|x)
$$
  and 
$$ p(y|x,  [\pmb{\varphi_{-3}}(x),a_3]) =
\int_{A_4(x)}\int_{A_2(x)}\int_{A_1(x)}p(y| x,a_1,a_2,a_3,a_4)
\varphi_1(da_1|x)\varphi_2(da_2|x)\varphi_4(da_4|x).$$
\\

\noindent{\bf Definition 2.} ({\bf Optimisation problem}) 
Let $ \pmb{\varphi}= (\varphi_1,\varphi_2,...,\varphi_n)\in \Phi.$
For each player $i\in {\cal N}$ consider the following  {\it constrained optimisation problem}
$COP(\pmb{\varphi_{-i}}):$
\begin{eqnarray}  \label{min} \nonumber
\lefteqn{\mbox{Minimise}}\\
&&\sum_{x\in X}\int_{A_i(x)}c_i^0(x,[\pmb{\varphi_{-i}}(x),a_i])
\mu(\{x\}\times da_i)
\end{eqnarray}
subject to $\mu\in \Pr(\mathbb{K}_i)$ and
\begin{equation}
\label{cop1}
\sum_{x\in X}\int_{A_i(x)}c_i^\ell(x,[\pmb{\varphi_{-i}}(x),a_i])
\mu(\{x\}\times da_i)\le \kappa^\ell_i \ \mbox{for all}\ \ell \in L, 
\end{equation}
and
\begin{equation}
\label{cop2}
\widehat{\mu}(x)= \mu(\{x\}\times A_i(x))=  (1-\alpha)\eta(x)+\alpha
\sum_{z\in X}\int_{A_i(z)}p(x|z,[\pmb{\varphi_{-i}}(z),a_i])\mu(\{z\}\times da_i) 
\end{equation}
 for all  $ x\in X.$ \\

We denote the set of all solutions to the problem 
$COP(\pmb{\varphi_{-i}})$ by    ${\cal O}_i(\pmb{\varphi_{-i}}).$\\

\noindent{\bf Remark 7.} (a)  In the  {\it constrained optimisation problem} $COP(\pmb{\varphi_{-i}})$  
player $i$   acts as the decision maker in a constrained 
discounted Markov decision  model. The transition probability, cost function and constraint functions are as follows: 
$p(x|z,[\pmb{\varphi_{-i}}(z),a_i]),$  $c_i^0(x,[\pmb{\varphi_{-i}}(x),a_i])$ 
and $c_i^\ell(x,[\pmb{\varphi_{-i}}(x),a_i]),$ $\ell\in L,$ respectively. Here,  $a_i\in A_i(x)$ and $x,\ z\in X.$
Equation (\ref{cop2}) implies that $\mu$ is an occupation measure
defined for this Markov decision process.
For details, see  \cite{bor}, Lemma 25 in \cite{piu} or Remark 6.3.1 in \cite{hll}.
Assumption {\bf C} assures that the set of all  occupation measures in $COP(\pmb{\varphi_{-i}})$  
satisfying  (\ref{cop1}) is non-empty.

(b) 
It is well-known from the literature  that  $COP(\pmb{\varphi_{-i}})$ 
has a solution in the bounded case under consideration. 
The set  ${\cal O}_i(\pmb{\varphi_{-i}}) $ is convex and compact. Indeed, 
the set of all    occupation measures, i.e., the measures satisfying  (\ref{cop2})  is convex and compact
(see for instance Theorem 3.1 in \cite{bor}).  Moreover,
since  $\mu\to\int_{\mathbb{K}_i}c_i^\ell(x,[\pmb{\varphi_{-i}}(x),a_i])
\mu(\{x\}\times da_i) $ is continuous on $\Pr(\mathbb{K}_i)$ for  all  $\ \ell \in L_0$, it follows that 
the subset of    occupation measures for which (\ref{cop1}) holds is closed,  and consequently compact.  
Hence, there exists an occupation measure that minimises  (\ref{min})  subject to (\ref{cop1}) and (\ref{cop2}). 
Thus, ${\cal O}_i(\pmb{\varphi_{-i}})$  is non-empty and  compact.
The convexity of ${\cal O}_i(\pmb{\varphi_{-i}})$  is obvious. \\

\noindent{\bf Remark 8.}   It should be noted that   
 \begin{equation}
\label{oimi} {\cal O}_i(\pmb{\varphi_{-i}})\subset  {\cal M}_i.
\end{equation} 
  If $\mu\in {\cal O}_i(\pmb{\varphi_{-i}})$, then, by 
 Lemma 3, there exists $\phi_i\in\Phi_i$ such that $\mu=\widehat{\mu}\phi_i$ where $\widehat{\mu}$ is the marginal of $\mu$ on $X.$
Furthermore,  $\mu(\{x\}\times\cdot)$ is the projection of 
$\rho \in {\cal M}$ defined as follows 
$$\rho(\{x\}\times d\pmb{a}):=\psi_1(da_1|x)\cdot\ldots\cdot\psi_n(da_n|x)\widehat{\mu}(x),$$ 
where $\psi_i=\phi_i$ and $\psi_j = \varphi_j$ for all $j\not= i.$ 
To see that $\rho$ is indeed an occupation measure from the set $\cal M$, it is sufficient to note from the definition of $\rho$ and      
 conclude from (\ref{cop2}) that 
\begin{eqnarray*}
& &\widehat{\rho}(x) = \widehat{\mu}(x)\\
&=&
(1-\alpha)\eta(x)+\alpha\sum_{z\in X}
\int_{A_n(z)}\ldots \int_{A_1(z)} p(x|z,a_1,\ldots,a_n)\psi_1(da_1|z)\ldots\psi_n(da_n|z) \widehat{\mu}(z)\\
&=&(1-\alpha)\eta(x)+\alpha \sum_{z\in X}
\int_{A_i(x)}p(x|z,[\pmb{\varphi_{-i}}(z),a_i])\mu(\{z\}\times da_i)\\
&=&(1-\alpha)\eta(x)+\alpha\sum_{z\in X}\int_{A(z)}p(x|z,\pmb{a})\rho(\{z\}\times d\pmb{a}).
\end{eqnarray*}
The claim now follows by applying \cite{bor}, or  Lemma 25 in \cite{piu},  or Remark 6.3.1 in \cite{hll}
\\
 
\noindent{\bf Remark 9.}  We would like to emphasise that the facts described in this remark hold
for games with arbitrary initial state distribution, unbounded functions $c_i^\ell$  
and transition probability satisfying our assumptions {\bf A} and {\bf B}.
For the details consult \cite{auto} and the literature mentioned below.
Assumption {\bf C} implies that for every multi-strategy $\pmb{\varphi}\in \Phi$ and each player $i\in {\cal N},$
there exists possibly non-stationary strategy  $\pi_i\in \Delta_i(\pmb{\varphi_{-i}}).$ 
Let $\mu_{\pi_i}$ be the  occupation measure defined as follows
$$
\mu_{\pi_i}(K) =  (1-\alpha)
\mathbb{E}_\eta^{\pi_i}\left[\sum_{t=1}^\infty \alpha^{t-1}  1_K(x^t,a_i^t)\right] \quad\mbox{
for each }K\in{\cal B}(\mathbb{K}_i).$$
Here, $ \mathbb{ E}_\eta^{\pi_i}$ denotes the expectation operator taken with respect to the unique probability 
measure defined  on the history space of the Markov decision  process 
governed by the transition probability $p(x|z,[\pmb{\varphi_{-i}}(z),a_i]),$ the initial distribution $\eta,$ 
and a strategy $\pi_i$ of the decision maker (player $i$). 
(Note that $\mathbb{ E}_\eta^{\pi_i} =\mathbb{E}_\eta^{[\pmb{\varphi_{-i}},\pi_i]}.$) 
From Proposition D.8 in \cite{hll}, it follows that  there exists a strategy $\sigma_i\in\Phi_i$ such that
$\mu_{\pi_i} =\widehat{\mu}_{\pi_i}\sigma_i,$ where  $\widehat{\mu}_{\pi_i}$ is the projection of $\mu_{\pi_i}$ on $X.$
Moreover,  from Lemma 3.1 in \cite{bor},  it follows that
$\mu_{\pi_i} =\mu_{\sigma_i} ,$ where $\mu_{\sigma_i} $ is an occupation measure defined as  above with
$\mathbb{ E}_\eta^{\pi_i}$ replaced by $\mathbb{ E}_\eta^{\sigma_i}.$
Therefore, by assumption {\bf C}, for every multi-strategy $\pmb{\varphi}\in \Phi$ and each player $i\in {\cal N},$
 there exists  $\sigma_i\in \Phi_i$ such that
\begin{eqnarray}\label{ch}
J_i^\ell([\pmb{\varphi_{-i}},\sigma_i])&=& \sum_{x\in X}\int_{A_i(x)}c_i^\ell(x,[\pmb{\varphi_{-i}}(x),a_i])
\mu_{\sigma_i}(\{x\}\times da_i)\\&=&\sum_{x\in X}\int_{A_i(x)}c_i^\ell(x,[\pmb{\varphi_{-i}}(x),a_i])
\mu_{\pi_i}(\{x\}\times da_i)=J_i^\ell([\pmb{\varphi_{-i}},\pi_i]) \nonumber
\end{eqnarray}
for all $\ell\in L.$ 
This equation implies that, if
$\pi_i \in   \Delta_i(\pmb{\varphi_{-i}}),$ then there exists a stationary strategy $\sigma_i\in \Phi_i$   such that 
$\sigma_i \in   \Delta_i(\pmb{\varphi_{-i}}),$ and we have
$$
J_i^\ell([\pmb{\varphi_{-i}},\pi_i]) =J_i^\ell([\pmb{\varphi_{-i}},\sigma_i]) < \kappa^\ell_i
$$ 
for all $\ell\in L.$ Since $\mu_{\pi_i}=\mu_{\sigma_i},$  the sequence of equalities in (\ref{ch})  is also valid for $\ell =0.$ \\

From now on,  we shall denote an element of ${\cal O}_i(\pmb{\varphi_{-i}})$  by $\mu_i.$  
When $\mu_i\in{\cal O}_i(\pmb{\varphi_{-i}})$, we take $\phi_i \in \Phi_i$ such that
$\mu_i= \widehat{\mu}_i\phi_i.$ Then, $\phi_i$ is an optimal stationary strategy for player $i$ 
in the constrained Markov decision  process   associated with   $COP(\pmb{\varphi_{-i}}).$ 
\\

\noindent{\bf Definition 3.} Under assumption {\bf D}, define the correspondence  $S:\prod_{i=1}^n {\cal M}_i
\to \prod_{i=1}^n {\cal M}_i$ by
$$S(\mu_1,\mu_2,...,\mu_n)= \prod_{i=1}^n {\cal O}_i(\pmb{\varphi_{-i}}),$$
where $\varphi_j\in \Phi_j$ is the unique strategy for player $j\in {\cal N}$
such that $\mu_j=\widehat{\mu}_j\varphi_j.$\\

By (\ref{oimi}) the correspondence $S$ is well-defined. 
We equip $\prod_{i=1}^n {\cal M}_i$ with the product topology.

The next result was proved in Lemma 2.1 in \cite{f}.\\

\noindent{\bf Lemma 5.} {\it Let {\bf A} and  {\bf B'} hold. 
The function $J_i^{\ell}(\cdot)$ is continuous on $\Phi$ for every  $\ell\in L_0.$ }\\

\noindent{\bf  Lemma 6.} {\it Assume {\bf A}, {\bf B'},  {\bf C} and {\bf D}. 
The correspondence $S$ is non-empty compact   convex-valued and is 
 upper semicontinuous.}\\

{\it Proof.} Since the spaces ${\cal M}_i$ are compact, to show the upper semicontinuity of the correspondence $S,$ 
it is enough to prove that $S$ has a closed graph. The other properties are discussed in Remark 7.  Assume that 
$\pmb{\mu}^k\to \pmb{\mu}^0$ in $\prod_{i=1}^n {\cal M}_i,$ where 
$\pmb{\mu}^k=(\mu_1^k,\ldots,\mu_n^k)$ for every $k\in\mathbb{N}\cup\{0\}.$
 Then, by Lemma 3, there exists  a unique strategy $\varphi_i^k\in \Phi_i$ such that
$\mu^k_i =\widehat{\mu}^k_i\varphi_i^k$ for every $i\in {\cal N}$ and  $k\in\mathbb{N}\cup\{0\}.$ 
Let $\nu^k_i\in {\cal O}_i(\pmb{\varphi}^k_{\pmb{-i}}),$
$k\in\mathbb{N},$ $i\in {\cal N}.$ 
Suppose that $\nu^k_i\to \nu^0_i$ weakly for every  $i\in {\cal N}.$ By Lemma 3, for any $i\in {\cal N}$ and $k\in\mathbb{N}\cup\{0\},$ 
there exists a unique strategy $\phi_i^k\in \Phi_i$ such that $\nu^k_i=\widehat{\nu}^k_i\phi_i^k,$ 
where $\widehat{\nu}^k_i$ is the marginal of 
 $\nu^k_i$ on $X.$  We have to show that  
\begin{equation}
\label{inclusion}
\nu^0_i\in {\cal O}_i(\pmb{\varphi}^0_{\pmb{-i}})\quad \mbox{ for every}\quad  i\in{\cal N}.
\end{equation}
From 
$\nu^k_i\in {\cal O}_i(\pmb{\varphi}^k_{\pmb{-i}}),$
for all $k\in\mathbb{N},$ $i\in {\cal N},$ it follows (see Remark 9) that
$$J_i^{\ell}([\pmb{\varphi}^k_{\pmb{-i}},\phi^k_i]) \le \kappa_i^\ell \quad\mbox{for all}\quad \ell\in L,\ i \in {\cal N}.$$
By Lemma 4, we know  that $\varphi_i^k\to \varphi_i^0$  and $\phi_i^k\to \phi_i^0$ in $\Phi_i$ for every $i\in{\cal  N}.$
 Thus, by Lemma 5, we conclude that 
\begin{equation}\label{a}
J_i^{\ell}([\pmb{\varphi}^0_{\pmb{-i}},\phi^0_i]) \le \kappa_i^\ell\quad \mbox{for all}\quad \ell\in L, \  i\in {\cal N}. \end{equation}
Moreover, we   have 
\begin{equation}\label{a0}
\lim_{k\to\infty}J_i^{0}([\pmb{\varphi}^k_{\pmb{-i}},\phi^k_i])=J_i^{0}([\pmb{\varphi}^0_{\pmb{-i}},\phi^0_i])\quad \mbox{for all}\quad i\in {\cal N}. 
\end{equation}
Inequality (\ref{a}) proves that the correspondence $\pmb{\varphi}_{\pmb{-i}}\to \Delta_i(\pmb{\varphi}_{\pmb{-i}})\cap \Phi_i$ 
has a closed graph.  Since all spaces $\Phi_j$ are compact, this correspondence   
is upper semicontinuous. 
By Lemma A2 in the Appendix, we conclude that $\pmb{\varphi}_{\pmb{-i}}\to \Delta_i(\pmb{\varphi}_{\pmb{-i}})\cap \Phi_i$ is continuous.
By the Berge maximum theorem, see pp. 115-116 in \cite{berge}, the function 
$$\pmb{\varphi}_{\pmb{-i}} \to   \min_{\sigma_i\in \Delta_i(\pmb{\varphi}_{\pmb{-i}})\cap \Phi_i}
J_i^{0}([\pmb{\varphi}_{\pmb{-i}},\sigma_i])$$
is continuous for any $i\in {\cal N}.$ Hence, 
\begin{equation}\label{a1}J_i^{0}([\pmb{\varphi}^k_{\pmb{-i}},\phi^k_i])=  \min_{\sigma_i\in \Delta_i(\pmb{\varphi}^k_{\pmb{-i}})\cap \Phi_i}
J_i^{0}([\pmb{\varphi}^k_{\pmb{-i}},\sigma_i])\ \to \   \min_{\sigma_i\in \Delta_i(\pmb{\varphi}^0_{\pmb{-i}})\cap \Phi_i}
J_i^{0}([\pmb{\varphi}^0_{\pmb{-i}},\sigma_i]) \quad\mbox{as}\ k\to\infty.\end{equation}
Expressions (\ref{a0}) and (\ref{a1}) 
imply that
$$
J_i^{0}([\pmb{\varphi}^0_{\pmb{-i}},\phi^0_i])
=
\min_{\sigma_i\in \Delta_i(\pmb{\varphi}^0_{\pmb{-i}})\cap \Phi_i} J_i^{0}([\pmb{\varphi}^0_{\pmb{-i}},\sigma_i])
\quad \mbox{for all}\quad i\in {\cal N}.
$$
This equality can be expressed in terms of problem $COP(\pmb{\varphi}^0_{\pmb{-i}}).$ Hence, (\ref{inclusion}) follows. 
 \hfill $\Box$\\

\noindent{\bf Proposition 1.} {\it Assume {\bf A}, {\bf B'},  {\bf C} and {\bf D}. Then   the $CSG$  
possesses a Nash equilibrium in  $\Phi.$}\\

{\it Proof.}  The set ${\cal M}_i$ can be viewed as a compact and convex subset of the  set, 
denoted by $ \mathbb{M}_i,$  of all signed finite measures on $X\times A_i$
equipped with the weak topology.  
$\mathbb{M}_i$ is a  locally convex topological Hausdorff 
space. Hence, the  set $\prod_{i=1}^n \mathbb{M}_i$ 
endowed with the product topology is also a locally convex topological  Hausdorff 
space. From Lemma  5 and the Kakutani-Fan-Glicksberg theorem (see Corollary 17.55 in \cite{ab}),  
it follows that there exists
$$(\mu_1^*,\mu_2^*,...,\mu_n^*)\in  S(\mu_1^*,\mu_2^*,...,\mu_n^*).$$
Now using  Lemma 3, take  $\varphi^*_i \in \Phi_i$ such that $\mu^*_i =
\widehat{\mu}^*_i\varphi_i^*$ for all $i\in {\cal N}.$ We claim that 
$\pmb{\varphi}^*= (\varphi^*_1, \ldots, \varphi^*_n)$ 
is a Nash equilibrium  in the
{\it CSG}.  We immediately have
\begin{equation}
\label{contr0}
J_i^{0}(\pmb{\varphi}^*) \le 
J_i^{0}([\pmb{\varphi}^*_{\pmb{-i}},\sigma_i])  \quad\mbox{for all}\quad
\sigma_i\in  \Phi_i \cap \Delta_i(\pmb{\varphi}^*_{\pmb{-i}}).
\end{equation}
Suppose that there exists some  $\pi_i\in \Pi_i$ such that
\begin{equation}
\label{contr1} J_i^{0}([\pmb{\varphi}^*_{\pmb{-i}},\pi_i])<  
J_i^{0}(\pmb{\varphi}^*) \quad\mbox{and}\quad\pi_i\in   \Delta_i(\pmb{\varphi}^*_{\pmb{-i}}).
\end{equation} 
Then, by Remark 9, we  conclude the existence of $\sigma_i\in\Phi_i$ for which
\begin{equation}
\label{contr2}
J_i^{0}([\pmb{\varphi}^*_{\pmb{-i}},\sigma_i])=J_i^{0}([\pmb{\varphi}^*_{\pmb{-i}},\pi_i]) 
\quad\mbox{and}\quad\sigma_i\in   \Delta_i(\pmb{\varphi}^*_{\pmb{-i}}).
\end{equation}
From (\ref{contr1}) and (\ref{contr2}) we get
$$
J_i^{0}([\pmb{\varphi}^*_{\pmb{-i}},\sigma_i]) < 
J_i^{0}(\pmb{\varphi}^*),
$$
which contradics (\ref{contr0}).
\hfill $\Box$

\section{ The proof of Theorem 1}

In this section, we introduce an  approximation  of the general game by  ones with truncated   cost and constraint functions and  
slightly perturbed initial state distributions. We apply Proposition 1 and other results from Section 4 to obtain 
Nash equilibria in the  truncated $CSGs$ and show that their limit is a Nash equilibrium in the original game.

Let $X=X_0\cup X_1,$ where $\eta(x)=0$ for $x\in X_0$ and $\eta(x)>0$ for $x\in X_1.$ 
Let $\widetilde{\eta}$ be a probability measure on $X_0$ such that $\widetilde{\eta}(x)>0$ for every $x\in X_0$.
For any $ m\in\mathbb{N}$,  define a perturbed initial state distribution  on $X$ as follows
\begin{equation}\label{p1}
\eta(m)(x):=\left(1-\frac{1}{m}\right)\eta(x) +\frac{1}{m} \widetilde{\eta}(x).\end{equation}
Clearly, $\eta(m)(x)>0$ and $\eta(m)(x)\to \eta(x)$   for every $x\in X$ as $m\to\infty.$

Moreover, for any $(x,\pmb{a})\in\mathbb{K},$ $i\in {\cal N},$  $\ell\in L_0$ and $m\in\mathbb{N},$   we set
\begin{equation}\label{p2}
c_i^{\ell,m}(x,\pmb{a}):=\left\{\begin{array}{ll}
-\sqrt{m},&\mbox{if } c_i^\ell(x,\pmb{a})<-\sqrt{m}\\
c_i^\ell(x,\pmb{a}),&\mbox{if } |c_i^\ell(x,\pmb{a})|\le \sqrt{m} \\
\sqrt{m},&\mbox{if } c_i^\ell(x,\pmb{a})>\sqrt{m}.\\
\end{array}\right.
\end{equation}

Before proving the theorem, we define a few functionals used in the proof. 
Let a multi-strategy  $\pmb{\pi}\in\Pi$ be fixed. For every $i\in {\cal N},$ $\ell\in L_0$ and $m\in\mathbb{N},$  put
\begin{eqnarray}\label{f1}
J_i^{\ell,m}(\pmb{\pi})
&:=& (1-\alpha)\mathbb{E}^{\pmb{\pi}}_{\eta}\left[\sum_{t=1}^\infty  \alpha^{t-1} c_i^{\ell,m}(x^t,\pmb{a}^t) \right],\\
\label{f2}
J_i^{\ell,\eta(m)}(\pmb{\pi})
&:=& (1-\alpha)\mathbb{E}^{\pmb{\pi}}_{\eta(m)}\left[\sum_{t=1}^\infty  \alpha^{t-1} c_i^{\ell,m}(x^t,\pmb{a}^t) \right].
\end{eqnarray}\\
Note that in (\ref{f1}) the initial state distribution is $\eta,$ while in (\ref{f2}) we use its perturbation $\eta(m).$ 
Both $J_i^{\ell,m}(\pmb{\pi})$ and $ J_i^{\ell,\eta(m)}(\pmb{\pi})$ 
are defined with the aid of truncated functions $c_i^{\ell,m}.$ 
 
The objective of  player $i\in{\cal N}$ in the modified game is to minimise 
$J_i^{0,\eta(m)}(\pmb{\pi})$ over $\pi_i\in\Pi_i$ with respect to the following constraints
$$J_i^{\ell,\eta(m)} (\pmb{\pi})\le  \kappa_i^{\ell,m},
$$ where  
\begin{equation}
\label{percon}
 \kappa_i^{\ell,m} :=
\left(1-\frac{1}{m}\right)\kappa_i^\ell+\frac{1}{\sqrt{ m}}
\quad\mbox{for all}\quad \ell\in L, \ i\in{\cal N}, \ m\in\mathbb{N}.
\end{equation} 

Denote by
$ J_i^{\ell,m}(x,\pmb{\varphi})$   the functional defined in (\ref{f1}) with the initial distribution $\eta$ 
replaced by the Dirac delta $\delta_{x}.$
Observe that  assumption {\bf C} holds with  $\kappa_i^{\ell,m}$ 
instead of $\kappa_i^\ell,$ since for any multi-strategy $\pmb{\varphi}\in\Phi$ there exists $\pi_i\in \Pi_i$ such that
\begin{eqnarray*}
J_i^{\ell,\eta(m)} ([\pmb{\varphi_{-i}},\pi_i])&=&
\left(1-\frac{1}{m}\right)\sum_{x\in X_1} J_i^{\ell,m} (x,[\pmb{\varphi_{-i}},\pi_i])\eta(x)+
\frac{1}{m}\sum_{x\in X_0} J_i^{\ell,m} (x,[\pmb{\varphi_{-i}},\pi_i])\widetilde{\eta}(x)\\ &<& 
\left(1-\frac{1}{m}\right)\kappa_i^\ell+\frac{1}{m}\sqrt{m}= \kappa_i^{\ell,m}
\end{eqnarray*}
for all $\ell\in L.$  

For any $\pmb{\pi}\in\Pi$ and $m\in\mathbb{N},$  define 
$$\Delta_i^{m}(\pmb{\pi_{-i}}):=\left\{\pi_i\in\Pi_i:\ J_i^{\ell,\eta(m)} (\pmb{\pi})\le
\kappa_i^{\ell,m}\quad\mbox{for all}\  \ell\in L\right\}.$$ \\

\noindent{\bf Definition 4.}  
The constrained discounted stochastic game with the  initial distribution (\ref{p1}), 
the cost and constraint functions as in (\ref{p2}),
the cost functionals as in  (\ref{f2}) and constants as in (\ref{percon}) is called an  $m$-$CSG.$ \\

\noindent{\bf Lemma 7.} {\it Let {\bf A} and  {\bf B} hold.
Then, the following holds.\\
(a) For every $\ell\in L_0$ and $i\in{\cal N}$
$$ \sup_{\pmb{\varphi}\in\Phi} \left|J_i^{\ell,\eta(m)}(\pmb{\varphi})- 
J_i^{\ell}(\pmb{\varphi})\right|\to 0\quad\mbox{as}\quad m\to\infty.$$
(b)  $ J_i^\ell(\cdot)$ 
is continuous on $\Phi$ for every  $\ell\in L_0$ and $i\in{\cal N}.$}\\

{\it Proof.}  (a) From the triangle inequality we have 
$$\sup_{\pmb{\varphi}\in\Phi} \left|J_i^{\ell,\eta(m)}(\pmb{\varphi})- 
J_i^{\ell}(\pmb{\varphi})\right|\le
 \sup_{\pmb{\varphi}\in\Phi} \left|J_i^{\ell, \eta(m)}(\pmb{\varphi})- 
J_i^{\ell,m}(\pmb{\varphi})\right|+ 
\sup_{\pmb{\varphi}\in\Phi} \left|J_i^{\ell,m}(\pmb{\varphi})- 
J_i^{\ell}(\pmb{\varphi})\right|=:\mbox{I + II}.$$
Then, we obtain 
\begin{eqnarray*}
\mbox{I }&=& \sup_{\pmb{\varphi}\in\Phi} \left|J_i^{\ell, \eta(m)}(\pmb{\varphi})- 
J_i^{\ell,m}(\pmb{\varphi})\right|
\le\sup_{\pmb{\varphi}\in\Phi}\sup_{y\in X}  
\left|J_i^{\ell,m}(y,\pmb{\varphi})\right|\cdot\sum_{x\in X}
|\eta(m)(x)-\eta(x)|\\
&\le& \sqrt{m}\left(  \sum_{x\in X_0} \frac{ \widetilde{\eta}(x)}{m} +   \sum_{x\in X_1} \frac{\eta(x)}{m} \right)=
 \frac{2}{\sqrt{m}}\to 0\quad\mbox{as } m\to\infty.
\end{eqnarray*}
Now let us consider the second term. 
By assumption {\bf B(i)} and (\ref{Bii1}), for any $\varepsilon>0,$  there exists $N_1\in\mathbb{N}$ such that for all $k> N_1,$
it holds
\begin{equation}\label{gw1}
\sup_{\pmb{\varphi}\in\Phi}\left|  (1-\alpha)\mathbb{E}^{\pmb{\varphi}}_{\eta}\left[\sum_{t=k}^\infty  \alpha^{t-1} c_i^{\ell}(x^t,\pmb{a}^t) \right]\right|
\le  \sup_{\pmb{\varphi}\in\Phi} (1-\alpha)\mathbb{E}^{\pmb{\varphi}}_{\eta}\left[\sum_{t=k}^\infty  \alpha^{t-1} w(x^t) \right] \le \frac{\varepsilon}4.
\end{equation}
Similarly, for all $m\in\mathbb{N},$ $k>N_1,$ $\ell\in L_0$ and $i\in{\cal N}$  
\begin{equation}\label{gw2}
\sup_{\pmb{\varphi}\in\Phi}\left|  (1-\alpha)\mathbb{E}^{\pmb{\varphi}}_{\eta}\left[\sum_{t=k}^\infty  \alpha^{t-1} 
c_i^{\ell,m}(x^t,\pmb{a}^t) \right]\right| \le \frac{\varepsilon}4.
\end{equation}
Note that
\begin{equation}\label{gw3}\left|  c_i^{\ell}(x,\pmb{a})- c_i^{\ell,m}(x,\pmb{a})   \right|\le w(x)1_{[ |c_i^\ell(x,\pmb{a})|> \sqrt{m}]}(x,\pmb{a}) \le 
w(x)1_{[ w(x)> \sqrt{m}]}(x).
\end{equation}
Consequently, by (\ref{gw1})-({\ref{gw3}), we obtain
\begin{eqnarray*}
\mbox{II }&=& \sup_{\pmb{\varphi}\in\Phi} \left|
J_i^{\ell,m}(\pmb{\varphi})- J_i^{\ell}(\pmb{\varphi})\right| \\
&\le& \sup_{\pmb{\varphi}\in\Phi}
\left|  (1-\alpha)\mathbb{E}^{\pmb{\varphi}}_{\eta}\left[\sum_{t=1}^{N_1}  \alpha^{t-1} 
c_i^{\ell,m}(x^t,\pmb{a}^t) \right]
- (1-\alpha)\mathbb{E}^{\pmb{\varphi}}_{\eta}\left[\sum_{t=1}^{N_1}  \alpha^{t-1} c_i^{\ell}(x^t,\pmb{a}^t)\right]\right|
+\frac{\varepsilon}2\\
&\le&  \sup_{\pmb{\varphi}\in\Phi}   (1-\alpha)\mathbb{E}^{\pmb{\varphi}}_{\eta}\left[\sum_{t=1}^{N_1}  \alpha^{t-1} 
w(x^t)1_{[ w(x^t)>\sqrt{m}]}\right]+\frac{\varepsilon}2
\end{eqnarray*}
for every $\ell\in L_0$ and $i\in{\cal N}.$
Hence,  from (\ref{Bii2}) for sufficiently large values of $m\in\mathbb{N},$ it follows  that
$$ \sup_{\pmb{\varphi}\in\Phi} \mathbb{E}^{\pmb{\varphi}}_{\eta}\left[w(x^t)1_{[ w(x^t)>\sqrt{m}]}\right]\le  \frac{\varepsilon}2$$
for all $t=1,\ldots, N_1.$ Therefore,
$$ \sup_{\pmb{\varphi}\in\Phi}   (1-\alpha)\mathbb{E}^{\pmb{\varphi}}_{\eta}\left[\sum_{t=1}^{N_1}  \alpha^{t-1} 
w(x^t)1_{[ w(x^t)>\sqrt{m}]}\right]\le  (1-\alpha)\sum_{t=1}^{N_1}  \alpha^{t-1}  \frac{\varepsilon}2 < \frac{\varepsilon}2.$$
Hence, $II<\varepsilon.$
This finishes the proof of part (a).

(b) By  Lemma 5, the functional  $J_i^{\ell,\eta(m)}(\cdot)$  is continuous on $\Phi$ for every $m\in\mathbb{N}$, $\ell\in L_0$ and 
$i\in {\cal N}.$ This fact and the uniform convergence proved in point (a)  imply the assertion. 
\hfill $\Box$\\

{\it Proof of Theorem 1.} From Proposition 1,  it follows that each 
$m$-$CSG$ possesses a   Nash equilibrium $\pmb{\phi}^{m}\in\Phi.$
Let  $(\pmb{\phi}^{m})_{m\in\mathbb{N}}$  be a   
sequence  in   $\Phi$ of  Nash equilibria in the $m$-$CSGs.$ From the compactness of $\Phi$, without loss of generality, we may assume that
$\pmb{\phi}^{m}$ converges to some $\pmb{\phi}\in \Phi$
as $m\to\infty.$ 
We claim that  $\pmb{\phi}$ is a Nash equilibrium in the original $CSG.$
Since $$ J_i^{\ell,\eta(m)}(\pmb{\phi}^{m})\le\kappa_i^{\ell,m}\quad\mbox{ for all } m\in\mathbb{N}, \ \ell\in L, \ i\in{\cal N},$$
by Lemma 7  and the fact that $\kappa_i^{\ell,m} \to\kappa_i^\ell$ as $m\to\infty$, it follows for  $\ell\in L$ and $i\in{\cal N}$ that
$$
 J_i^{\ell,\eta(m)}(\pmb{\phi}^{m})\to J_i^\ell(\pmb{\phi}) \le \kappa_i^\ell
\quad \mbox{ as}\quad m\to\infty.$$
These facts  immediately entail that $\pmb{\phi}\in\Phi$ is feasible in the original $CSG.$ The rest will follow, if we show
that
$$  J_i^0(\pmb{\phi})=\min_{\pi_i\in\Delta_i(\pmb{\phi_{-i}})}J_i^0([\pmb{\phi_{-i}},\pi_i])
\quad\mbox{for every}\  i\in{\cal N}.$$
On the contrary, assume that there exists   player $i\in{\cal N}$ and a strategy $\bar{\pi}_i\in \Pi_i$ such that $\bar{\pi}_i\in  \Delta_i(\pmb{\phi_{-i}})$ and
$$J_i^0([\pmb{\phi_{-i}}, \bar{\pi}_i])< J_i^0(\pmb{\phi}).$$
By Remark 9, we may replace $\bar{\pi}_i$ by  a strategy $\gamma_i\in  \Phi_i$  in the sense that  $\gamma_i\in  \Delta_i(\pmb{\phi_{-i}})$ and
\begin{equation}\label{a2} 
J_i^0([\pmb{\phi_{-i}}, \bar{\pi}_i])=J_i^0([\pmb{\phi_{-i}}, \gamma_i])< J_i^0(\pmb{\phi}).\end{equation}
By Lemma A1 in the Appendix, one can select a sequence 
$(\gamma_i^m)_{m\in\mathbb{N}}$ with $\gamma_i^m\in \Delta_i^{m}
(\pmb{\phi}^{m}_{\pmb{-i}})$ such that
\begin{equation}\label{q}
\lim_{m\to\infty}J_i^{0,\eta(m)} ([\pmb{\phi}^{m}_{\pmb{-i}},\gamma^m_i]) =  J_i^{0}([\pmb{\phi_{\pmb{-i}}},\gamma_i]). 
\end{equation}
Since  $\phi^{m}_i$ is the $i$-th coordinate of the Nash equilibrium profile $\pmb{\phi}^{m}$ in the $m$-$CSG$, we have
$$J_i^{0,\eta(m)}([\pmb{\phi}^{m}_{\pmb{-i}},\phi^{m}_i])=   \min_{\sigma_i\in \Delta_i^{m}(\pmb{\phi}^{m}_{\pmb{-i}})}  
J_i^{0,\eta(m)}([\pmb{\phi}^{m}_{\pmb{-i}},\sigma_i])\le J_i^{0, \eta(m)}([\pmb{\phi}^{m}_{\pmb{-i}},\gamma^m_i]).$$
Taking the limit as $m\to\infty$ in the above display and applying Lemma 7 and (\ref{q}), we get
$$J_i^0(\pmb{\phi})= J_i^{0}([\pmb{\phi}_{\pmb{-i}},\phi_i])\le J_i^{0}([\pmb{\phi}_{\pmb{-i}},\gamma_i]).$$
This inequality  contradicts (\ref{a2}). 
 \hfill $\Box$\\

In the above proof we tacitly assumed that  $X_0\not= \emptyset.$  If  $X_0 = \emptyset,$  
then our   proof can be simplified in an obvious manner.

\section{Additional  remarks on assumptions}

In this section, we give some examples and comments on assumptions {\bf B} and {\bf W}. 
For simplicity, we   consider a one-person game, i.e., a constrained discounted Markov decision process,
where the player is called a decision maker. Therefore, 
$A_1(x)=A(x)$ for all $x\in X$   and an  element of $A(x)$  will be denoted by $a$ instead of $\pmb{a}.$

In the following  example, inspired by the example of Blackwell  \cite{bl},  
the function $w$ satisfying assumption {\bf W(ii)}   does not exist. \\

\noindent{\bf  Example 1.}  We consider a simple Markov decision process.
Let $X=\mathbb{N}\cup \mathbb{N}^*,$ where $\mathbb{N}^*:=\{1^*,2^*,3^*,\ldots\}.$ 
The action sets are: $A(n)=\{c,s\}$  and $A(n^*)=\{s\}$
for $n\in\mathbb{N}$, $n^*\in\mathbb{N}^*.$ Here, $c$ means {\it continue} and $s$ means {\it stop.} 
State $1^*$ is absorbing and $p(1^*|n^*,s)=p((n+1)^*|n,s)=1$ for $n\in\mathbb{N},$ $n^*\in\mathbb{N}^*.$   
Moreover,
$p(1^*|n,c)=1-p(n+1|n,c)=q,$ $n\in\mathbb{N}$, 
$n^*\in\mathbb{N}^*,$ where $q\in [0,1].$ 
The cost functions are non-negative and satisfy inequalities:
$c(1^*,s)=0\le c(n,s)=c(n,c)\le 1$ and $c(n^*,s)\le n$ for $n\in\mathbb{N},$  $n^*\in\mathbb{N}^*.$
 
We begin with showing that our assumption {\bf B} holds. We define the function $w$ in the simplest way, i.e., 
$$w(n)=1,\quad w(n^*)=n, \quad\mbox{for }\ n\in\mathbb{N},\  n^*\in\mathbb{N}^*.$$
The initial distribution is geometric and is given on the set $\mathbb{N}$:  $\eta(m)=(1-g)g^{m-1}$ for each $m\in \mathbb{N},$  
where $g\in (0,1) $ is fixed. 
We now prove {\bf B(ii)}. Fix any strategy $\pi$ of the decision maker. Denote by $\mathbb{E}^\pi_m$  the expectation operator 
on the trajectories of the process governed by $\pi$ and starting at state $m\in \mathbb{N}.$  
We note that
\begin{eqnarray*}
\mathbb{E}^\pi_\eta\left(\sum_{k=n}^\infty \alpha^{k-1}w(x^k)\right)&=& 
\sum_{m=1}^\infty \left[  \mathbb{E}^\pi_m\left(\sum_{k=n}^\infty \alpha^{k-1}w(x^k)\right)  \right]
(1 -g)g^{m-1}\\
&\le&  \sum_{m=1}^\infty \left[  \mathbb{E}^\pi_m\left(\sum_{k=n}^\infty \alpha^{k-1}(k+m-1)\right)  \right]
(1-g)g^{m-1}\\
&=& \frac{\alpha^{n-1}(n-1)(1-\alpha)+\alpha^n}{(1-\alpha)^2}+\frac{\alpha^{n-1}}{g(1-\alpha)}\to 0\quad\mbox{as } n\to\infty.
\end{eqnarray*}
Hence, (\ref{Bii1}) is satisfied. The inequality in the above display is due to the observation that  $w(x^k)\le k+(m-1)$ if the initial state
is $m$ and  no matter which strategy the decision maker  uses. Now we show  that (\ref{Bii2})  holds. 
Fix $t\in\mathbb{N}$ in  (\ref{Bii2}) and any strategy $\pi$. Then,
$$
\mathbb{E}^\pi_\eta\left( w(x^t)1_{[w(x^t)\ge k]} \right)\le  \sum_{m=1}^\infty \left[ (t+m-1) \mathbb{E}^\pi_m 1_{[w(x^t)\ge k]} \right]g^{m-1}(1-g).$$
For any $\varepsilon>0,$ there exists $M\in\mathbb{N}$ such that
 $$\max\left\{   (t-1)\sum_{m=M}^\infty g^{m-1}(1-g),  \sum_{m=M} ^\infty mg^{m-1}(1-g) \right\} \le \frac\varepsilon 2. $$
Set $k:=t+M-1$ and note that
$$  \mathbb{E}^\pi_m 1_{[w(x^t)\ge k]}  \le   1_{[t+m-1\ge k]}=   1_{[m\ge M]}.$$
Hence, for any strategy $\pi$ of the decision maker, we have
$$
\mathbb{E}^\pi_\eta\left( w(x^t)1_{[w(x^t)\ge k]} \right)\le   \sum_{m=M}^\infty (t+m-1)g^{m-1}(1-g)\le \frac \varepsilon 2+ \frac \varepsilon 2=\varepsilon.$$
This proves  (\ref{Bii2}).

Now it  can be easily seen that the inequality in {\bf W(i)} does not hold. Indeed,  there is no $\delta>1$ such that the inequality 
 $$\sum_{y\in X} w(y)p(y|n,s) =w((n+1)^*)=n+1 \le \delta w(n) =\delta$$ holds for all
$n\in\mathbb{N}.$ 

The second possibility is to change a function $w$ in such a  way that $w(n)=n$. 
Then, we obtain the special conditions on the discount factor $\alpha$ through  inequalities from {\bf W(i)}:
$$\sum_{y\in X} w(y)p(y|n,s) =w((n+1)^*)=n+1 \le n\delta \ \Rightarrow \ 
1+\frac 1n\le\delta  \ \Rightarrow \ \delta=2\mbox{ and }\alpha<\frac 12, $$
and 
$$ \sum_{y\in X} w(y)p(y|n,c)=q+(1-q)(n+1)\le n \delta  \ \Rightarrow \ 1-q +\frac 1n\le \delta 
\ \Rightarrow \ \delta=2-q\mbox{ and }\alpha<\frac 1{2-q}. $$
Thus, {\bf W} holds, but only for $\alpha<1/2.$ 
This is a serious restriction for the discount factor. 
Other inequalities in {\bf W(i)} are automatically satisfied and we do not consider them here.
Finally, the third possibility is to modify $w$ by adding some constant $d>0$ (but we keep $w(1^*)=1$).  
This idea was first discussed in \cite{jn}.
The above inequalities are as follows:
$$\sum_{y\in X} w(y)p(y|n,s) =w((n+1)^*)=n+1+d \le (n+d)\delta  \ \Rightarrow \ 1+\frac 1{n+d}\le\delta$$
and  $$ \sum_{y\in X} w(y)p(y|n,c)=q+(1-q)(n+1+d)\le (n+d) \delta  \ \Rightarrow \ 1-q+\frac 1{n+d}\le\delta$$
for all $n\in\mathbb{N}.$ Hence, for any $\alpha\in (0,1),$ 
we may choose a constant $d$ such that $\alpha\delta<1$ with $\delta:=1+\frac 1{1+d}.$
This forces us to select carefully an appropriate function $w$. The second disadvantage is that the function $w$
is less natural than the original $w$ given  above, i.e., when $w(n)=1$ for all $n\in \mathbb{N}.$  \hfill $\Box$

\section{Appendix}

First note that from (\ref{percon}), it follows that for sufficiently large values of $m$ we have  
\begin{equation}
\label{kapm}
\kappa_i^\ell \le \kappa_i^{\ell,m}\quad\mbox{for all}\quad \ell \in L, \ i\in {\cal N}.
\end{equation}
Indeed, this inequality is equivalent to $\kappa_i^\ell \le 
\sqrt{m}.$

The following result  is crucial in the proofs of Lemma 6 and Theorem 1.\\ 
 \\
\noindent{\bf Lemma A1.} {\it Let {\bf A}, {\bf B} and  {\bf C} 
hold  and let $\pmb{\phi}\in \Phi$ be feasible in the $CSG.$ Assume that $\pmb{\phi}^m\to 
 \pmb{\phi}$ in $\Phi$ as $m\to\infty.$ 
  For any  $i\in{\cal N}$, 
$\gamma_i\in \Delta_i(\pmb{\phi_{-i}})$,   
 there exists a sequence of strategies $(\gamma_i^{m})_{m\in\mathbb{N}}$ in $\Phi_i$ such that $\gamma_i^{m}\in  
\Delta_i^{m}(\pmb{\phi}^{m}_{\pmb{-i}})$  for every $m\in\mathbb{N}$
and \begin{equation}
\label{limj0}
J_i^{0,\eta(m)}
([\pmb{\phi}^m_{\pmb{-i}},\gamma^m_i]) \to  J_i^{0}
([\pmb{\phi}_{\pmb{-i}},\gamma_i])\quad\mbox{as}\quad     m\to\infty.
\end{equation}}
\\

{\it Proof.} Recall that, if $\pmb{\phi}^m = (\phi^m_1,...,\phi_n^m)$ and
$\pmb{\phi}  = (\phi_1,...,\phi_n), $ then $\pmb{\phi}^m \to \pmb{\phi}$ means that
$\phi_j^m  \to \phi_j$ for all $j\in {\cal N}$ as $m\to\infty.$  

Observe that from Lemma 7, it follows that, for every $\ell\in L$, $i\in {\cal N},$
\begin{equation}
\label{limj1}
J_i^\ell([\pmb{\phi_{-i}},\gamma_i])=\lim_{m\to\infty} J_i^{\ell,\eta(m)}
([\pmb{\phi}^m_{\pmb{-i}},\gamma_i])\le\kappa_i^\ell.
\end{equation}
We consider two cases. 

 1. If  $J_i^\ell([\pmb{\phi_{-i}},\gamma_i])<\kappa_i^\ell$ for all $\ell\in L$, 
then by (\ref{kapm}) there exists $N\in\mathbb{N}$ such that, 
for all $m>N$ and for all $\ell\in L,$ 
$$J_i^{\ell,\eta(m)}
([\pmb{\phi}^m_{\pmb{-i}},\gamma_i])\le\kappa_i^\ell \le \kappa_i^{\ell,m}. $$
Hence, $\gamma_i\in  \Delta_i^{m}(\pmb{\phi}^m_{\pmb{-i}})$ for all $m>N.$ Therefore, it suffices to put
$\gamma_i^m:=\gamma_i$ for $m>N.$ For $m=1,\ldots, N$ from assumption {\bf C},  it is enough to
take any $\gamma_i^m\in  \Delta_i^{m}(\pmb{\phi}^m_{\pmb{-i}})\cap \Phi_i.$
Then (\ref{limj0}) follows immediately from (\ref{limj1}). 

2. Let $i\in {\cal N}$ be fixed. Assume now that for at least one $\ell\in L$,  
$J_i^\ell([\pmb{\phi_{-i}},\gamma_i])=\kappa_i^\ell.$ By assumption
{\bf C} and Remark 9, there exists a strategy $\xi_i\in \Phi_i$ for which 
$$J_i^\ell([\pmb{\phi_{-i}}, \xi_i])< \kappa_i^\ell\quad\mbox{for all}\quad \ell\in L.$$ 
Hence, by Lemma 7, there is a constant $\kappa>0$ such that 
\begin{equation}\label{a21}
J_i^\ell([\pmb{\phi_{-i}}, \xi_i])=\lim_{m\to\infty}  
J_i^{\ell,\eta(m)}
([\pmb{\phi}^m_{\pmb{-i}},\xi_i])\le  \kappa_i^\ell- \kappa\end{equation}
and 
\begin{equation}\label{a22} \lim_{m\to\infty}  
J_i^{\ell,\eta(m)}
([\pmb{\phi}^m_{\pmb{-i}},\gamma_i])\le  \kappa_i^\ell
\end{equation}
for all $\ell\in L. $    

Let  $\mu^m_{\gamma_i}$  be the  occupation measure defined on $\mathbb{K}_i,$
when the Markov process is induced by the initial distribution $\eta(m),$ 
the transition probability $p(x|z,[\pmb{\phi}^m_{\pmb{-i}}(z),a_i])$ 
and the strategy $\gamma_i$ of the decision maker (player $i$).
By definition
$$
\mu^m_{\gamma_i}(K) =  (1-\alpha)
\mathbb{E}_{\eta(m)}^{[\pmb{\phi}^m_{\pmb{-i}},\gamma_i]}\left[\sum_{t=1}^\infty \alpha^{t-1} 
1_K(x^t,a_i^t)\right]\quad\mbox{
for any }K\in{\cal B}(\mathbb{K}_i). 
$$
Analogously, we define $\mu^m_{\xi_i}.$
Thus, from Remark 9, (\ref{a21}) and (\ref{a22}), we deduce that,
for every $\tau\in (0,\kappa),$ there exists $N_\tau\in\mathbb{N}$ such that
\begin{equation}\label{wz1}
\sum_{x\in X}\int_{A_i(x)} c_i^{\ell,m}(x,[\pmb{\phi}^m_{\pmb{-i}}(x),a_i])\mu^m_{\xi_i}(\{x\}\times da_i)<
\kappa_i^\ell-\kappa+\tau
 \end{equation}
 and
\begin{equation}\label{wz2}
\sum_{x\in X}\int_{A_i(x)} c_i^{\ell,m}(x,[\pmb{\phi}^m_{\pmb{-i}}(x),a_i])\mu^m_{\gamma_i}(\{x\}\times da_i)
<\kappa_i^\ell+\tau
\end{equation}
 and
\begin{equation}\label{wz3}
0< \kappa_i^{\ell,m}-\kappa_i^{\ell}<\tau
\end{equation}
for all $m>N_\tau$ and $\ell\in L.$ 

Put
$$\lambda(\tau,m):=\frac{\tau+\kappa_i^\ell-\kappa_i^{\ell,m}}{\kappa}$$
and note by (\ref{wz3}) that $\lambda(\tau,m)\in(0,1)$ and $\lambda(\tau,m)  \searrow 0$ as $\tau\searrow 0.$
Moreover,  for all $m>N_\tau$ and $\ell\in L$, it follows from (\ref{wz1}) and (\ref{wz2}) that
\begin{eqnarray}\label{fs1}
\lefteqn{\left(1-\lambda(\tau,m)\right) \sum_{x\in X}\int_{A_i(x)} c_i^{\ell,m}(x,[\pmb{\phi}^m_{\pmb{-i}}(x),a_i])
\mu^m_{\gamma_i}(\{x\}\times da_i) } \\&&+\lambda(\tau,m)\sum_{x\in X}\int_{A_i(x)} 
c_i^{\ell,m}(x,[\pmb{\phi}^m_{\pmb{-i}}(x),a_i])
\mu^m_{\xi_i}(\{x\}\times da_i) \nonumber  \\ 
 &<& 
 \left(1-\lambda(\tau,m)\right)(\kappa_i^\ell+\tau)+\lambda(\tau,m)(\kappa_i^\ell-\kappa+\tau)  
= \kappa_i^\ell +\tau-\lambda(\tau,m)\cdot\kappa  
 =  \kappa_i^{\ell,m}.  \nonumber 
\end{eqnarray}

Let $(\epsilon_k)_{k\in\mathbb{N}}$ be a sequence of numbers in $(0,\kappa)$ such that $\epsilon_k \searrow 0$ as $k\to \infty.$  
For each $ \tau=\epsilon_k,$ there exists $N_{k}$ such that 
(\ref{wz1})-(\ref{wz3}) hold for all $m> N_{k}.$  
We may assume that the sequence $(N_{k})_{k\in\mathbb{N}}$ is 
increasing, so $\lim_{k\to\infty}N_{k}=\infty.$  Note that, for each $m> N_1,$ there exists a unique $k$ such that $N_{k} < m
\le N_{k+1}.$ Using this positive integer $k,$  we put $\tau(m)=\epsilon_k.$ 
Observe that $\tau(m) = \epsilon_k$ for all $m$ such that $N_{k} < m
\le N_{k+1}.$

Let
$$\lambda(\tau(m),m):=\frac{\tau(m)+\kappa_i^\ell-\kappa_i^{\ell,m}}{\kappa}.$$
Note that, if $m\to\infty,$ then $N_{k}    \to
\infty$
and $\tau(m)\to 0.$  Thus, $\lambda(\tau(m),m) \to 0$ as $m\to\infty.$ 
 
Now we define a new occupation measure as follows: 
for each $m$ such that $ N_{k}<m\le N_{k+1},$ we set
\begin{equation}\label{fs2}
\nu_i^m=\left(1-\lambda(\tau(m),m)\right) \mu^m_{\gamma_i}+\lambda(\tau(m),m) \mu^m_{\xi_i}, \end{equation}
where  $\tau(m)=\epsilon_k.$ 
Observe that $\mu^m_{\gamma_i} $   and $\mu^m_{\xi_i} $  belong to ${\cal M}_i.$ 
By Lemma 2, $\nu_i^m \in {\cal M}_i.$  
Moreover,  (\ref{fs1}) holds for
$\tau =\tau(m)=\epsilon_k$ and for all $m$ such that $ N_{k}<m\le N_{k+1}.$
By Lemma 3, for each $m\in\mathbb{N},$ there exists a unique strategy $\gamma^m_i\in\Phi_i$  
such that   $\nu_i^m=\widehat{\nu}_i^m\gamma_i^m.$
Note that by (\ref{fs1}), (\ref{fs2})  and Remark 9, we have 
$$
J_i^{\ell,\eta(m)}([\pmb{\phi}^{m}_{\pmb{-i}},\gamma_i^m])\le \kappa_i^{\ell,m}
$$
for all $\ell\in L$ and   for all $m$ such that $N_{k}  <m\le 
N_{k+1}.$ 
Hence, $\gamma_i^m\in \Delta_i^{m}(\pmb{\phi}^{m}_{\pmb{-i}})\cap \Phi_i$ for all  $m>N_{ 1}.$ Clearly,
for $m=1,\ldots, N_{1}$ assumption {\bf C} enables us to choose  
any $\gamma_i^m\in  \Delta_i^{m}(\pmb{\phi}^m_{\pmb{-i}})\cap \Phi_i.$

Let $\ell=0.$ Making use of Remark 9 and Lemma 7,  we infer  that
\begin{eqnarray}\label{j1}
\lefteqn{  J_i^{0}([\pmb{\phi}_{\pmb{-i}},\gamma_i]) =
\lim_{m\to\infty}J_i^{0,\eta(m)}([\pmb{\phi}^{m}_{\pmb{-i}},\gamma_i])}\\\nonumber
&&=\lim_{m\to\infty}
\sum_{x\in X}\int_{A_i(x)}c_i^{0,m}(x,[\pmb{\phi}^m_{\pmb{-i}}(x),a_i])
\mu^m_{\gamma_i}(\{x\}\times da_i).
\end{eqnarray}
Since $\lambda(\tau(m),m)\to 0$ in   
 (\ref{fs2}) as $m\to\infty$,  we have
\begin{eqnarray}\label{j2}
\lefteqn{ \lim_{m\to\infty}
\sum_{x\in X}\int_{A_i(x)}c_i^{0,m}(x,[\pmb{\phi}^m_{\pmb{-i}}(x),a_i])
\mu^m_{\gamma_i}(\{x\}\times da_i)}\\
&&=\lim_{m\to\infty}
\sum_{x\in X}\int_{A_i(x)}c_i^{0,m}(x,[\pmb{\phi}^m_{\pmb{-i}}(x),a_i])
\nu^{m}_{i}(\{x\}\times da_i). \nonumber
\end{eqnarray}
However, $\nu^{m}_{i}=\widehat{\nu}^m_i\gamma_i^m$ and  
$\nu^{m}_{i}$ is a convex combination of two occupation measures
$\mu^m_{\gamma_i}$  and   $\mu^m_{\xi_i}$ determined (among others) 
by the disturbed initial state distribution $\eta(m)$. 
Therefore,
\begin{equation}\label{j3}
\sum_{x\in X}\int_{A_i(x)}c_i^{0,m}(x,[\pmb{\phi}^m_{\pmb{-i}}(x),a_i])
\nu^{m}_{i}(\{x\}\times da_i)= J_i^{0,\eta(m)}([\pmb{\phi}^{m}_{\pmb{-i}},\gamma^m_i]).\end{equation}
Consequently, combing together (\ref{j1}), (\ref{j2}) and (\ref{j3}), we conclude that
$$\lim_{m\to\infty}J_i^{0,\eta(m)}([\pmb{\phi}^{m}_{\pmb{-i}},\gamma^m_i])=J_i^{0}([\pmb{\phi}_{\pmb{-i}},\gamma_i]).
$$
This finishes the proof.   \hfill $\Box$\\

The next result can be proved in a similar  manner as the above lemma
with some necessary amendments. Lemma A2 is used in Lemma 6,  where we assume that the 
cost and constraint functions are bounded and the support of the initial distribution is the whole state space $X$.\\

\noindent{\bf Lemma A2. }  {\it Let {\bf A}  and  {\bf C} hold.
Assume that for each $i\in {\cal N}$ and $\ell \in L$ the 
function $c^\ell_i$ is bounded and $\eta(x)>0$ for all $x\in X.$  
Then   the   correspondence $\pmb{\phi_{-i}} \to \Delta_i(\pmb{\phi_{-i}})\cap \Phi_i$ from 
$\Phi_{-i} =\prod_{j\not=i}\Phi_j$ to $\Phi_i$  
is lower semicontinuous 
for each player $i\in{\cal N}.$}\\

{\it Proof.} 
We have to prove that, if
$\pmb{\phi}_{\pmb{-i}}\in \Phi_{\pmb{-i}},$
$\gamma_i \in \Delta_i(\pmb{\phi_{-i}})\cap \Phi_i$
and  $\pmb{\phi}^{m}_{\pmb{-i}}\to \pmb{\phi}_{\pmb{-i}}$ in $\Phi_{\pmb{-i}}$
as $m\to\infty,$ then there exists a sequence 	
$(\gamma_i^{m})_{m\in\mathbb{N}}$ in $\Phi_i$ such that $\gamma_i^{m}\in  
\Delta_i(\pmb{\phi}^{m}_{\pmb{-i}})$  for every $m\in\mathbb{N}$ and 
$ \gamma^m_i \to \gamma_i$  as $m\to\infty.$

Observe that from Lemma 5, it follows that, for every $\ell\in L$, $i\in {\cal N},$
$$
\lim_{m\to\infty} J_i^{\ell}
([\pmb{\phi}^m_{\pmb{-i}},\gamma_i])= J_i^\ell([\pmb{\phi_{-i}},\gamma_i]) \le\kappa_i^\ell.
$$
We consider two cases. 

1. If  $J_i^\ell([\pmb{\phi_{-i}},\gamma_i])<\kappa_i^\ell$ for all $\ell\in L$, 
then there exists $N\in\mathbb{N}$ such that, 
for all $m>N$ and for all $\ell\in L,$ 
$$J_i^{\ell}
([\pmb{\phi}^m_{\pmb{-i}},\gamma_i])\le\kappa_i^\ell. $$
Hence, $\gamma_i\in  \Delta_i(\pmb{\phi}^m_{\pmb{-i}})$ for all $m>N.$ Therefore, it suffices to put
$\gamma_i^m:=\gamma_i$ for $m>N.$ For $m=1,\ldots, N$ from assumption {\bf C},  it is enough to
take any $\gamma_i^m\in  \Delta_i(\pmb{\phi}^m_{\pmb{-i}})\cap \Phi_i.$
Then, the convergence of $\gamma_i^m$ to $\gamma_i$  is obvious.

2. Fix $i\in {\cal N}.$  Assume that for at least one $\ell\in L$, it holds  that  
$J_i^\ell([\pmb{\phi_{-i}},\gamma_i])=\kappa_i^\ell.$ From assumption
{\bf C} and Remark 9, there exists a strategy $\xi_i\in \Phi_i$ for which 
$$J_i^\ell([\pmb{\phi_{-i}}, \xi_i])< \kappa_i^\ell\quad\mbox{for all}\quad \ell\in L.$$ 
Hence, by Lemma 5, there is a constant $\kappa>0$ such that 
$$
J_i^\ell([\pmb{\phi_{-i}}, \xi_i])=\lim_{m\to\infty}  
J_i^{\ell}
([\pmb{\phi}^m_{\pmb{-i}},\xi_i])\le  \kappa_i^\ell- \kappa
$$
and 
$$ \lim_{m\to\infty}  
J_i^{\ell}
([\pmb{\phi}^m_{\pmb{-i}},\gamma_i])\le  \kappa_i^\ell
$$
for all $\ell\in L. $    

Define  the occupation measures $\mu^m_{\gamma_i}$  and $\mu^m_{\xi_i}$ 
on $\mathbb{K}_i$ as in Case 2 in the proof of Lemma A1, but with $\eta(m)$ instead of $\eta.$ 
In addition, replace   $c_i^{\ell,m}$ by $ c_i^\ell$, $\kappa^{\ell,m}_i$ by $\kappa^{\ell}_i$ and 
$\lambda(\tau,m)$ by $\lambda(\tau) :=\tau/\kappa $  in the proof of Lemma A2. 
Then,  we introduce the occupation measures $\nu_i^m$ in a similar way as in (\ref{fs2}). Namely,
\begin{equation}
\label{nuim}
\nu_i^m=\left(1-\lambda(\tau(m))\right) \mu^m_{\gamma_i}+\lambda(\tau(m))
 \mu^m_{\xi_i},
\quad\mbox{if}\quad N_{k} <m\le N_{k+1},\quad  k=\tau(m).
\end{equation} 
The definition of $\tau(m)$ is the same as in the proof of Lemma A1, but 
$\lambda(\tau(m)) = \tau(m)/\kappa.$ 
From the proof of Lemma A1 (or Lemma 3),
we also conclude that, for any $m\in\mathbb{N},$  there exists a unique  $\gamma_i^m\in \Phi_i$ such that   
$\nu^m_i= \widehat{\nu}^m_i\gamma_i^m.$ 
Moreover, 
$\gamma_i^{m}\in  
\Delta_i(\pmb{\phi}^{m}_{\pmb{-i}})$  for every $m\in\mathbb{N}.$
It remains to show that $\gamma_i^m \to \gamma_i$ as $m\to\infty.$

Let  $\mu_{\gamma_i}$  be the  occupation measure defined on $\mathbb{K}_i,$   
when the Markov process is induced by the initial distribution $\eta,$ 
the transition probability $p(x|z,[\pmb{\phi}_{\pmb{-i}}(z),a_i])$ and   $\gamma_i.$ 
For any bounded continuous function $f\in C(\mathbb{K}_i),$  we put
$$\widehat{f}(x ) =\int_{A_i(x )}f(x,a_i )\gamma_i(da_i|x).$$
 Then, for each $f\in C(\mathbb{K}_i),$ we have 
\begin{eqnarray*}
\lefteqn{ \sum_{x\in X}\int_{A_i(x)}f(x,a_i)\mu_{\gamma_i}(\{x\}\times da_i) 
=  (1-\alpha)
 \mathbb{E}_{\eta }^{[\pmb{\phi}_{\pmb{-i}},\gamma_i]}\left[\sum_{t=1}^\infty \alpha^{t-1} 
 f(x^t,a_i^t)\right]}\\ &=&
(1-\alpha)
 \mathbb{E}_{\eta }^{[\pmb{\phi}_{\pmb{-i}},\gamma_i]}\left[\sum_{t=1}^\infty \alpha^{t-1} 
 \widehat{f}(x^t )\right]=
\sum_{x\in X}\int_{A_i(x)}\widehat{f}(x )\mu_{\gamma_i}(\{x\}\times da_i)\\ &=&
\sum_{x\in X} \widehat{f}(x )\widehat{\mu}_{\gamma_i}(x) =
\sum_{x\in X} \int_{A_i(x )}f(x ,a_i) 
\gamma_i(da_i|x)\widehat{\mu}_{\gamma_i}(x).
\end{eqnarray*}
Since $f\in C(\mathbb{K}_i)$ is arbitrary, it follows that  
$ 
\mu_{\gamma_i} =\widehat{\mu}_{\gamma_i}\gamma_i.
$ 
By Lemma  2.1 in \cite{f} (or Lemma 5 with $c^\ell_i=f$), for each $f\in C(\mathbb{K}_i),$ we have 
\begin{eqnarray*}
&&\sum_{x\in X}\int_{A_i(x)}f(x,a_i)\mu^m_{\gamma_i}(\{x\}\times da_i) 
=  (1-\alpha)
\mathbb{E}_{\eta }^{[\pmb{\phi}^m_{\pmb{-i}},\gamma_i]}\left[\sum_{t=1}^\infty \alpha^{t-1} 
f(x^t,a_i^t)\right] \\  &\to &
(1-\alpha)
\mathbb{E}_{\eta }^{[\pmb{\phi}_{\pmb{-i}},\gamma_i]}\left[\sum_{t=1}^\infty \alpha^{t-1} 
f(x^t,a^t_i )\right]=
\sum_{x\in X}\int_{A_i(x)}f(x,a_i )\mu_{\gamma_i}(\{x\}\times da_i).
\end{eqnarray*}
Thus, $\mu^m_{\gamma_i}$  converges weakly to $\mu_{\gamma_i}$ as $m\to\infty.$ 
This fact and (\ref{nuim}) imply that 
$\nu^m_{i}= \widehat{\nu}^m_{i}\gamma^m_i$ 
converges weakly to $\mu_{\gamma_i}= \widehat{\mu}_{\gamma_i}\gamma_i  $ as $m\to\infty.$ 
Since $\eta(x) >0$ for all $x\in X,$ by Lemma 4, $\gamma^m_i\to \gamma_i$ in $\Phi_i$ as $m\to\infty.$ 
This finishes the proof.   \hfill $\Box$\\

\end{document}